\documentclass[pre,twocolumn]{revtex4}

\usepackage{amsmath}
\usepackage{amsfonts}
\usepackage{hyperref}
\usepackage{graphicx}

\usepackage[T1]{fontenc}

\renewcommand{\L}{\mathrm{Li}}%
\renewcommand{\S}{\mathcal{S}}%
\renewcommand{\d}{\mathrm{d}}%

\begin{document}

\title{Solving Kepler's equation via nonlinear sequence transformations}

\author{Riccardo Borghi}
\affiliation{ Dipartimento di Ingegneria, Universit\`{a} ``Roma Tre'', Via Vito Volterra 62, I-00146 Rome, Italy}

\begin{abstract}
Since more than three centuries Kepler's equation continues to represents an important benchmark for testing new computational techniques.
In the present paper, the classical Kapteyn series solution of Kepler's equation originally conceived by Lagrange and Bessel will be revisited from a 
different  perspective, offered by the relatively new and still largely unexplored framework of the so-called nonlinear sequence 
transformations.
 
The main scope of the paper is to provide numerical evidences supporting the fact that the Kapteyn series solution of Kepler's equation could be a Stieltjes 
series. To support such a conjecture, two types of Levin-type sequence transformations, namely  Levin $d$- and  Weniger $\delta$-transformations, 
will be employed to sum up several wildly divergent series derived by the  Debye representation of Bessel functions. As an interesting byproduct of this analysis, an 
effective recursive algorithm to generate arbitrarily higher-order Debye's polynomials will be developed.

Such a  ``Stieltjeness'' conjecture will also be numerically validated by directly employing the Levin-type transformations to accelerate the complex 
Kapteyn series solution of the Kepler equation. Both $d$- and $\delta$- transformations display exponential convergence, whose rate will be 
numerically estimated.  A few conclusive words, together with some hints for future extensions in the direction of more general class of Kapteyn series, 
eventually close the paper.
 
\end{abstract}

\maketitle

\section{Introduction}
\label{Sec:Intro}

Likely, there are a very few equations in the world (if any) which can boast such a large 
number of solving strategies as much as the celebrated {Kepler equation} (KE henceforth).
Since 1650 a huge number of different methods have been conceived to numerically solve KE. 
In the classical textbook by Colwell~\cite{Colwell/1993}, most of such  methods 
have been listed and briefly described. However, since 1993 such a list dramatically 
increased and it is really hard to give account for all strategies developed so far. 
Just to give an idea, it is worth quoting  some of the papers dealing with KE which have been published within the last ten years~\cite%
 {Eisinberg2010215,Gu2010716,Davis201059,Antonov2011182,Liu2012,Yang2012517,%
Calvo2013143,Farnocchia201321,Mortari20141,Avendano201427,Badolati2015316,Lynden-Bell2015363,%
Avendano2015435,Oltrogge2015271,Rauh20161,Ebaid2016457,Calvo201719,Raposo-Pulido20171702,%
Chanda201,Boetzel/Susobhanan/Gopakumar/Klein/Jetzer/2017,Ebaid20171,Alshaery2017933,%
Alshaery201727,Elipe2017415,Lppez20182583,Aljohani2018,Orlando2018849,Raposo-Pulido2018,%
Zechmeister2018,Ibrahim20192269,Calvo2019,Tommasini2020,Abubekerov/Gostev/2020,Sacchetti/2020,Tommasini2020b,Zechmeister2021,Gonzalez/Hernandez/2021,Tommasini/2021,Philcox/Goodman/Slepian/2021}.

Despite its apparent mathematical simplicity, solving KE continues to play a pivotal role 
in the science of computation. As it was pointed out again in~\cite{Colwell/1993}, 
\begin{quotation}
{\em Any new technique for the treatment of trascendental equations should be applied to this illustrious case; any new 
insight, however slight, lets its conceiver join an eminent list of contributors.
}
\end{quotation}

In the present paper we intend to give a further contribution to the subject, by focusing our attention  on a celebrated semi-analytic  
approach to solve  KE,  based on a Bessel function series expansion~\cite{Colwell/1993}. 
A similar  approach has recently been applied to search solutions of the so-called post-Newtonian 
KE's~\cite{Boetzel/Susobhanan/Gopakumar/Klein/Jetzer/2017}.
Although the Bessel solution of KE has been abandoned for practical, computational purposes, 
it presents interesting unexplored features that, in our opinion, deserve to be investigated. 

The Lagrange/Bessel solution of KE was the first example of a class of expansions called \emph{Kapteyn series} (KS henceforth)~\cite{Kapteyn/1893},
which in recent times have gained a role of considerable importance~\cite{Eisinberg/2010,Tautz/Lerche/Dominici/2011,Baricz/Jankov/Pogany/2011,Nikishov/2014,Dragana/Masirevic/Pogany/2017,Xue/Li/Man/Xing/Liu/Li/Wu/2019}. 
One of the  scopes of the present paper is to provide numerical evidences suggesting that a particular type of  Kapteyn series, which is strictly 
connected to the Bessel  solution of KE,  is a \emph{Stieltjes series}. 
This, in turns, would imply the whole computational machinery developed for resumming 
Stieltjes asymptotic series to be employable in order to retrieve the KE solution. This will be done in the second part of the paper, where two examples of the so-called Levin-type nonlinear sequence 
transformation 
will be used to accomplish such a computational task.

The paper is structured as follows: in Sec.~\ref{Sec:TheoAnal}  a brief resume of the Bessel 
solution of KE's equation is outlined, together with a few words about Stieltjes functions and Stieltjes series.
Section~\ref{Sec:StieltjesKapteyn} is devoted to numerically explore the Stieltjes nature of the Bessel KS,
while in Sec.~\ref{Sec:KeplerWT} the results of several numerical experiments carried out on directly applying  Levin-type transformations on the 
Bessel KS are shown. Some conclusive words, together with some hints to extend our approach to deal with more general classes of Kapteyn series, 
are given in Sec.~\ref{Sec:Conclusions}.

\section{Preliminaries}
\label{Sec:TheoAnal}

\subsection{The Bessel solution of Kepler's equation}
\label{Sec:BesselKE}

\begin{figure}[!ht]
\includegraphics[width=6cm]{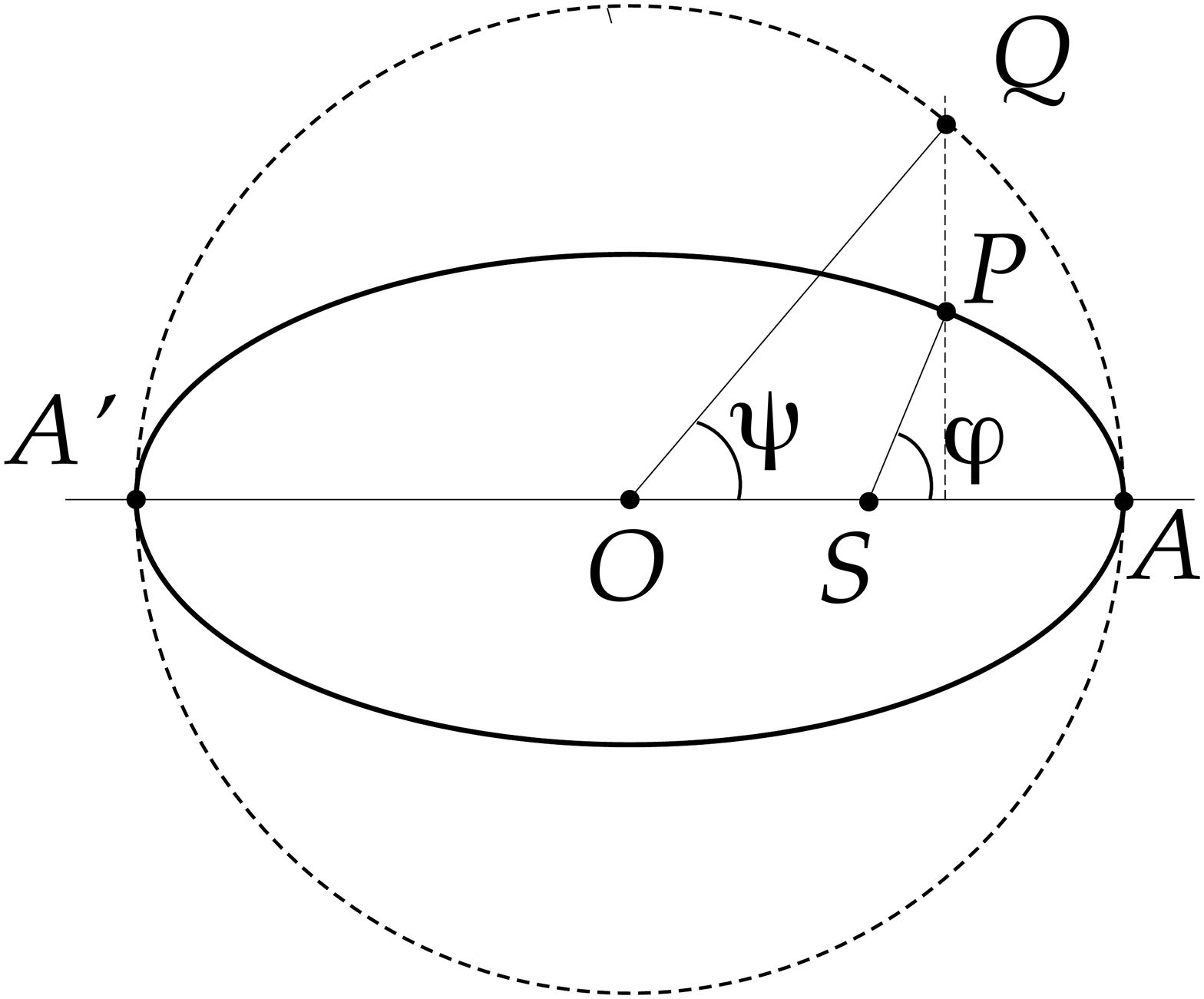}
\caption{The geometry of Kepler's Equation.}
\label{Fig:KEGeometry}
\end{figure}

The (elliptic) KE is formally defined as follows:
\begin{equation}
\label{Eq:KE.1}
M\,=\,\psi\,-\,\epsilon\,\sin\,\psi\,,
\end{equation}
where the meaning of the quantities $M$, $\psi$, and $\epsilon$ can be grasped on referring to Fig.~\ref{Fig:KEGeometry}.

Consider the motion of a planet around the star $S$ along an elliptic orbit with eccentricity $\epsilon$~\cite{Orlando2018849}.
Suppose $P$ to be the position of the planet at the time $t$, which is measured on assuming 
the apocenter $A$  as the initial (i.e., at $t=0$) position (motion is counterclockwise). The mean  angular speed $\omega=2\pi/T$, with $T$ 
being the orbital period, is  introduced in such a  way that $M=\omega t$, which is called the \emph{mean anomaly}.
To retrieve the position of the planet along the orbit at time $t>0$, the angle $\psi$ is first evaluated by solving KE~(\ref{Eq:KE.1}), which gives
 the position of the point $Q$ along the circumscribed circle. $P$ is then immediately obtained  by the geometric 
construction depicted in Fig.~\ref{Fig:KEGeometry}.

In the present paper our attention will be focused on the following  explicit Fourier series representation of the KE solution:~\cite[Ch.~3]{Colwell/1993}:
\begin{equation}
\label{Eq:KE.2}
\begin{array}{l}
\displaystyle
\psi\,=\,M\,+\,S(\epsilon;M)\,,
\end{array}
\end{equation}
where the function $S(\epsilon;M)$ is defined as
\begin{equation}
\label{Eq:KE.3}
\begin{array}{l}
\displaystyle
S(\epsilon;M)\,=\,\sum_{n=1}^\infty\,\frac{2\,J_n(n\,\epsilon)}n\,\sin\,nM\,.
\end{array}
\end{equation}
%
%
Although the  above series converges  for any  $\epsilon\in[0,1)$, it turns out that such a  convergence is extremely slow, especially when $\epsilon \to 1$. 
To give a single numerical example of such a slow convergence,  in Tab.~\ref{Tab:SlowConvergence} the behaviour of the partial sums of the series
into~Eq.~(\ref{Eq:KE.3}) are shown as a function of their order for $\epsilon=9/10$ and $M=\pi/4$.
%
\begin{table}[!ht]
\centerline{
    \begin{tabular}{|c|c|}
    \hline
    Partial sum order & Estimate of $\psi$\footnote{The exact value given by numerically solving Eq.~(\ref{Eq:KE.1})
    is 1.6800337357880455291\ldots} \\ \hline \hline
	0  	& 	1.35949\ldots \\ \hline
	1	&	1.66564\ldots \\ \hline
	2	&	1.78539\ldots \\ \hline
	3	&	1.78539\ldots \\ \hline
	4	&	1.73032\ldots \\ \hline
	5	&	1.67194\ldots \\ \hline
	10	&	1.70076\ldots  \\ \hline
	15	&	1.66772\ldots \\ \hline
	20	&	1.68367\ldots \\ \hline
	25	&	1.68138\ldots \\ \hline
	30	&	1.67725\ldots \\ \hline
	35	&	1.68210\ldots \\ \hline
	40	&	1.67933\ldots \\ \hline
	45	&	1.67968\ldots \\ \hline
	50	&	1.68076\ldots\\ \hline
	55	&	1.67945\ldots\\ \hline
	60	&	1.68023\ldots\\ \hline
	65	&	1.68014\ldots\\ \hline
	70	&	1.67978\ldots\\ \hline
    \hline
    \end{tabular}
}
\caption{Behaviour ot the partial sums of the Fourier series into~Eq.~(\ref{Eq:KE.3})  as a function of their order for $\epsilon=9/10$ and $M=\pi/4$.}
\label{Tab:SlowConvergence}
\end{table}

It is not difficult to understand why such bad convergence features, together with the fact that a considerable number of Bessel function evaluations has to be 
implemented, unavoidably brought the Fourier series to be definitely abandoned as far as practical applications of KE are concerned. 
Nevertheless,  the Bessel series expansion~(\ref{Eq:KE.3}) presents features that make it a subject of considerable 
interest in math as well as  in theoretical physics. The scope of the present paper is to explore some of these features.
To this end, function $S$ is first recast as the 
imaginary part of the complex  function $\S(\epsilon;M)$ defined through the following KS:
\begin{equation}
\label{Eq:KE.4}
\begin{array}{l}
\displaystyle
\S(\epsilon;M)\,=\,\sum_{n=1}^\infty\,\frac{2\,J_n(n\,\epsilon)}n\,\exp(\mathrm{i}nM)\,,
\end{array}
\end{equation}
which is nothing but a complex power series expansion evaluated across the unit circle.
Its  convergence could be proved, for example, on invoking the following  Bessel function asymptotics, valid for  
$\epsilon < 1$~\cite[Sec.~8.4]{Watson/1995}:
\begin{equation}
\label{Eq:KE.4.1}
\begin{array}{l}
\displaystyle
J_n(n\,\epsilon)\,\sim\,\frac 1{\sqrt{2\pi\,\sqrt{1\,-\,\epsilon^2}}}\,\frac{\rho^n}{\sqrt n}\,,\qquad\quad  n\to\infty\,,
\end{array}
\end{equation}
where
\begin{equation}
\label{Eq:KE.4.2}
\begin{array}{l}
\displaystyle
\rho\,=\,\exp\left(\sqrt{1\,-\,\epsilon^2}\right)\,\frac{1\,-\,\sqrt{1\,-\,\epsilon^2}}\epsilon\,,
\end{array}
\end{equation}
is a positive parameter less than unity when~$\epsilon\in(0,1)$.
On substituting from Eq.~(\ref{Eq:KE.4.1}) into Eq.~(\ref{Eq:KE.4}) and on introducing the complex parameter 
$z=\rho\,\exp(\mathrm{i}M)$, it turns out that the single terms of the series~(\ref{Eq:KE.4}) asymptotically approach 
those of the paradigmatic model series
\begin{equation}
\label{Eq:KE.6}
\begin{array}{l}
\displaystyle
\sum_{n=1}^\infty\,\frac{z^n}{n^{3/2}}\,,
\end{array}
\end{equation}
which, for $|z|<1$, is nothing but the Dirichlet series representation of so-called polylogarithm function $\L_{3/2}(z)$, where 
\begin{equation}
\label{Eq:KE.6.1}
\begin{array}{l}
\displaystyle
\L_{\nu}(z)\,=\,\sum_{n=1}^\infty\,\frac{z^n}{n^{\nu}}\,,\qquad |z| < 1\,.
\end{array}
\end{equation}
For $|z| \ge 1$, the series into Eq.~(\ref{Eq:KE.6.1}) can be analytically extended to the whole complex plane but the
half-axis $\mathrm{Re}\{z\} > 1$. More importantly, polylogarithm $\L_{\nu}(z)$ is an example of \emph{Stieltjes function}. 
For reader's convenience, the main definitions of Stieltjes functions and of Stieltjes series will now be briefly summarized. 
Interested readers can found some useful references for instance in Ref.~\cite{Caliceti/Meyer-Hermann/Ribeca/Surzhykov/Jentschura/2007}.

A function $f(z)$, with $z\in\mathbb{C}$,  is  called Stieltjes function if it admits the following integral representation:
\begin{equation}
\label{Eq:Stieltjes.1}
\displaystyle
{f}(z)=
\displaystyle\int_0^\infty\,\frac{\mathrm{d}\mu(t)}{t+z}\,,\qquad |\arg(z)|<\pi\,,
\end{equation}
where $\mu(t)$ is a bounded, \emph{nondecreasing} function taking infinitely
many different values on the interval $0 \le t < \infty$. 
The integral into Eq.~(\ref{Eq:Stieltjes.1}) is called \emph{Stieltjes integral}.
On introducing  the moment sequence $\{\mu_m\}^\infty_{m=0}$, whose elements are defined as
\begin{equation}
\label{Eq:Stieltjes.2}
  \mu_m \; = \; \int_{0}^{\infty} \, t^m \, \mathrm{d} \mu\,,
  \qquad m \ge 0 \, ,
\end{equation}
it turns out that they are nonnegative for all values of $m$. 
Any Stieltjes integral~(\ref{Eq:Stieltjes.1}) can be associated to a so-called \emph{Stieltjes series}, which is asymptotic to $f(z)$ as $z\to\infty$, in the 
Poincar\'e sense, as follows:
\begin{equation}
\label{Eq:Stieltjes.3}
\displaystyle
f(z) \sim \sum_{m=0}^\infty\,\frac {(-1)^m\,\mu_m}{z^{m+1}}\,,\qquad  z \to \infty\,.
\end{equation}
Whether such a series converges or diverges depends ultimately on the asymptotic behavior of the
moment sequence  $\{\mu_m\}^\infty_{m=0}$ as $m\to\infty$.

%
%
%
%

Now, proving  the  Stieltjes nature of polylogarithm functions is almost trivial.
It is sufficient to start from the following integral representation of $\L_\nu$~\cite[25.12.11]{DLMF}:
\begin{equation}
\label{Eq:Debye.1.1.2.2}
\begin{array}{l}
\displaystyle
\L_\nu(z)\,=\,\dfrac z{\Gamma(\nu)}\,
\int_0^\infty\,
\dfrac{x^{\nu-1}}{\exp(x)-z\,}\,\d x\,=\,\\
\\
\displaystyle
\,=\,
-\dfrac 1{\Gamma(\nu)}\,
\int_0^\infty\,
\dfrac{x^{\nu-1}\exp(-x)}{\exp(-x)-1/z\,}\,\d x,
\end{array}
\end{equation}
and set  $t=\exp(-x)$ to have 
\begin{equation}
\label{Eq:Debye.1.1.2.3}
\begin{array}{l}
\displaystyle
\L_\nu(z)\,=\,
-\dfrac 1{\Gamma(\nu)}\,
\int_0^1\,
\dfrac{\left(-\log t\right)^{\nu-1}}{t\,-\,1/z}\,
\mathrm{d}t\,.
\end{array}
\end{equation}
%
Apart from the overall minus sign, the function $\L_\nu$ can then be written in the form of an Stieltjes integral~(\ref{Eq:Stieltjes.1}), 
with $-1/z$ in place of $z$ and with the measure, say $\mu_\nu(t)$, defined as
\begin{equation}
\label{Eq:Debye.1.1.2.4}
\mu_\nu(t)\,=\,
\left\{
\begin{array}{lr}
\displaystyle
\dfrac{\Gamma(\nu,-\log t)}{\Gamma(\nu)}\,, & \qquad \qquad 0\le t \le 1,\\
&\\
1\,& \qquad \qquad t>1\,,
\end{array}
\right.
\end{equation}
where symbols $\Gamma(\cdot)$ and $\Gamma(\cdot,\cdot)$ denote the gamma and the  \emph{incomplete} gamma functions, respectively~\cite{DLMF}.

The above described ``asymptotic connection'' between series~(\ref{Eq:KE.4}) and~(\ref{Eq:KE.6}) via Eq.~(\ref{Eq:KE.4.1}) would 
at first sight suggest the function $\mathcal{S}(\epsilon;M)$ to be Stieltjes too.
Although we are not able to provide a rigorous proof about such a conjecture, in the next section some strong numerical evidences supporting it will be given.

\section{Is $\mathcal{S}(\epsilon;M)$ a {Stieltjes} series?}
\label{Sec:StieltjesKapteyn}

\subsection{Debye's expansion}
\label{Subsec:PreStieltjes}

The Bessel function asymptotics given into Eq.~(\ref{Eq:KE.4.1}) is the leading term of 
a complete series expansion known as {Debye's expansion}. More explicitly,  we have~\cite[formula~(10.19.3)]{DLMF}
\begin{equation}
\label{Eq:Debye.1}
\begin{array}{l}
\displaystyle
{J_n(n\,\epsilon)}\,=\,
\dfrac {\rho^n}{\sqrt{2\pi\,\sqrt{1\,-\,\epsilon^2}}}\,
\sum_{k=0}^\infty\,\dfrac 1{n^{k+1/2}}\,U_k\,\left(\dfrac 1{\sqrt{1-\epsilon^2}}\right)\,,
\end{array}
\end{equation}
where the symbol $U_k(\cdot)$ denotes the so-called $k$th-order Debye polynomial.

Debye polynomials are still largely unexplored mathematical objects. They are formally defined through the following 
integro-differential recurrence rule~\cite[formula~(10.41.9)]{DLMF}:
\begin{equation}
\label{Eq:Debye.1.1}
\begin{array}{l}
\displaystyle
U_{k+1}(t)\,=\,
\dfrac{t^2}2(1-t^2)\dfrac{\d U_k}{\d t}\,+\,\dfrac 18\,\int_0^t\,(1-5\xi^2)\,U_k(\xi)\,\d\xi\,,
\end{array}
\end{equation}
with the initial condition $U_0(t)=1$.
Now, from a mere computational point of view, definition~(\ref{Eq:Debye.1.1}) makes the generation of higher-order Debye 
polynomials considerably time demanding, even if powerful symbolic platforms are employed. 
To give an idea, the commercial software \emph{Mathematica 12.0} requires, on a DELL workstation,  about 10~seconds to generate $U_{11}(t)$ but
more than 12~minutes to generate $U_{13}(t)$!

Since a deep numerical analysis of the Debye expansion~(\ref{Eq:Debye.1}) is mandatory for the scopes of the present paper, 
in Appendix~\ref{App:B} an effective recursive algorithm to compute only the coefficients of $U_k(t)$, starting from the definition into
Eq.~(\ref{Eq:Debye.1.1}), has been outlined. In this way, to generate Debye polynomials having degrees of the order of hundreds or 
even of thousands takes only a few seconds. For instance, to produce $U_{1000}(t)$ on the above described computing 
machine and with the same symbolic platform only 25~seconds are required.

\begin{figure}[!ht]
\centering
\begin{minipage}[t]{7cm}
\includegraphics[width=6cm]{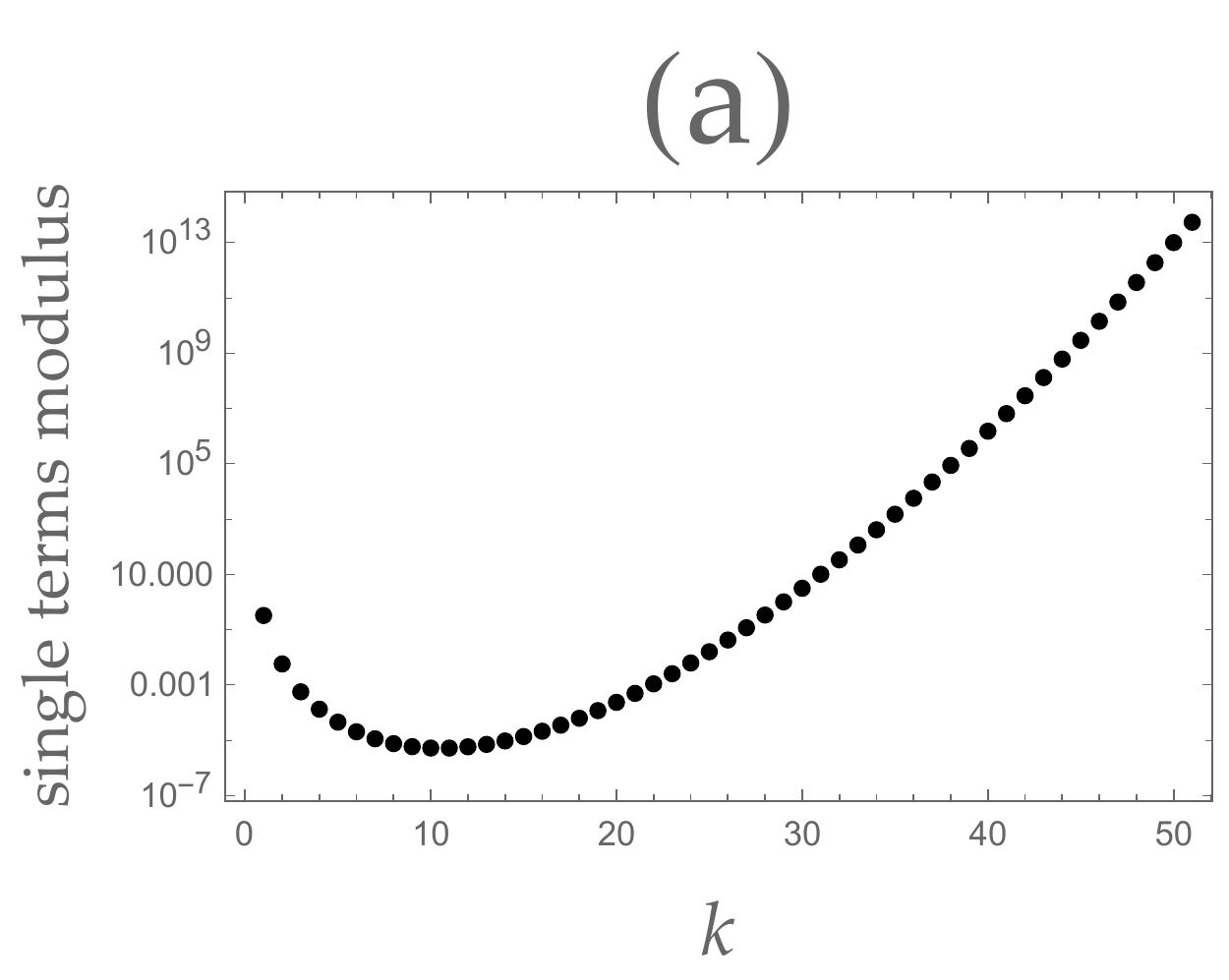}
\hfill
\includegraphics[width=6cm]{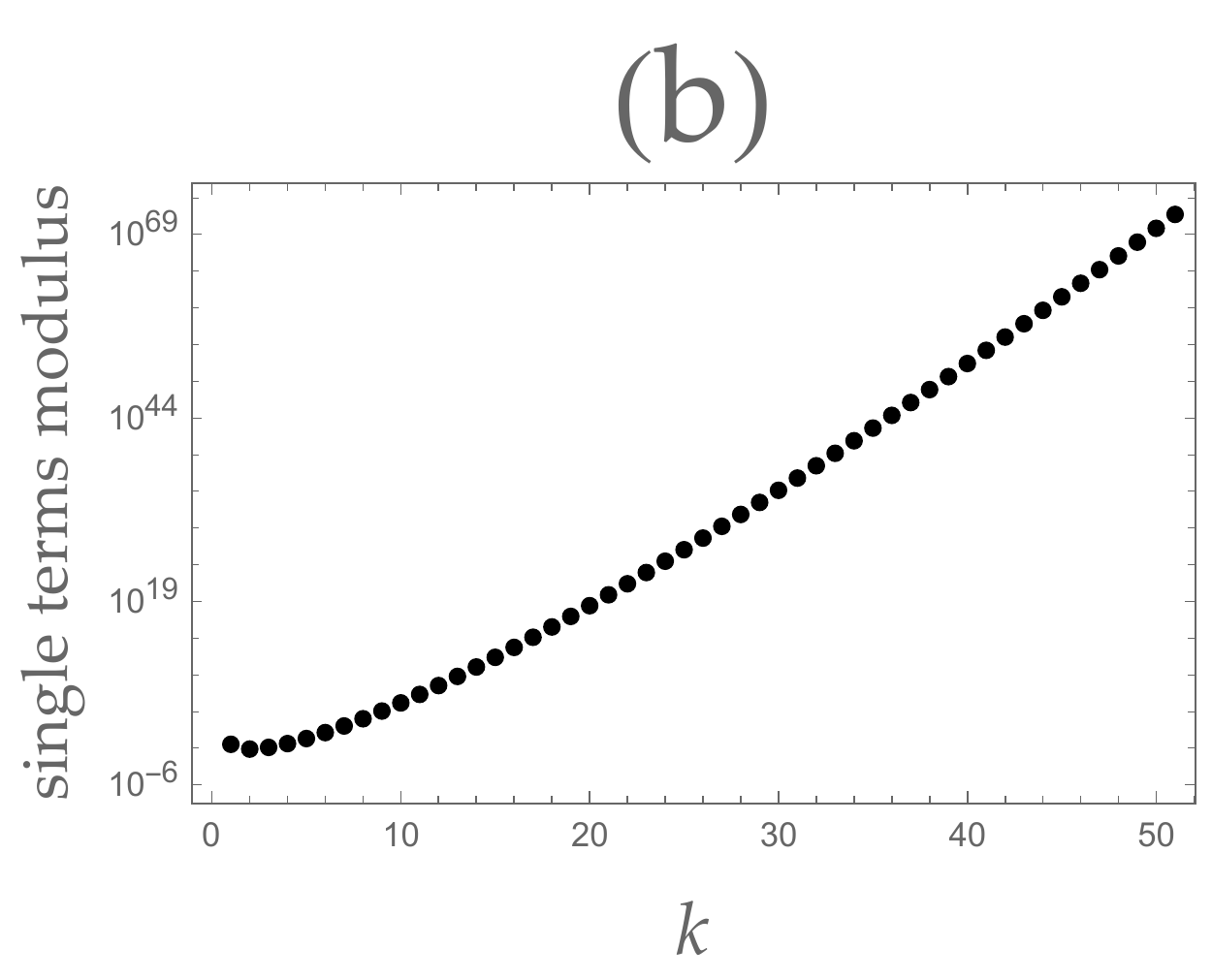}
\end{minipage}
\caption{Behavior, as a function of $k$,  of the modulus of the single terms of the series in 
Eq.~(\ref{Eq:Debye.1}), shown for $n=10$ and $\epsilon=5/10$ (a), $\epsilon=9/10$ (b).}
\label{Fig:Debye.1}
\end{figure}

Equipped with such a powerful computational tool, we are now ready to start our numerical exploration of
the Debye expansion~(\ref{Eq:Debye.1}). 
To begin, in Fig.~\ref{Fig:Debye.1} the behavior, as a function of $k$, of the modulus of the single terms of the series into Eq.~(\ref{Eq:Debye.1}) is 
shown for $n=10$ and for $\epsilon=5/10$ (a) and $\epsilon=9/10$ (b). 
Both sequences appear to be clearly divergent. 

\begin{figure}[!ht]
\centering
\begin{minipage}[t]{8cm}
\includegraphics[width=8cm]{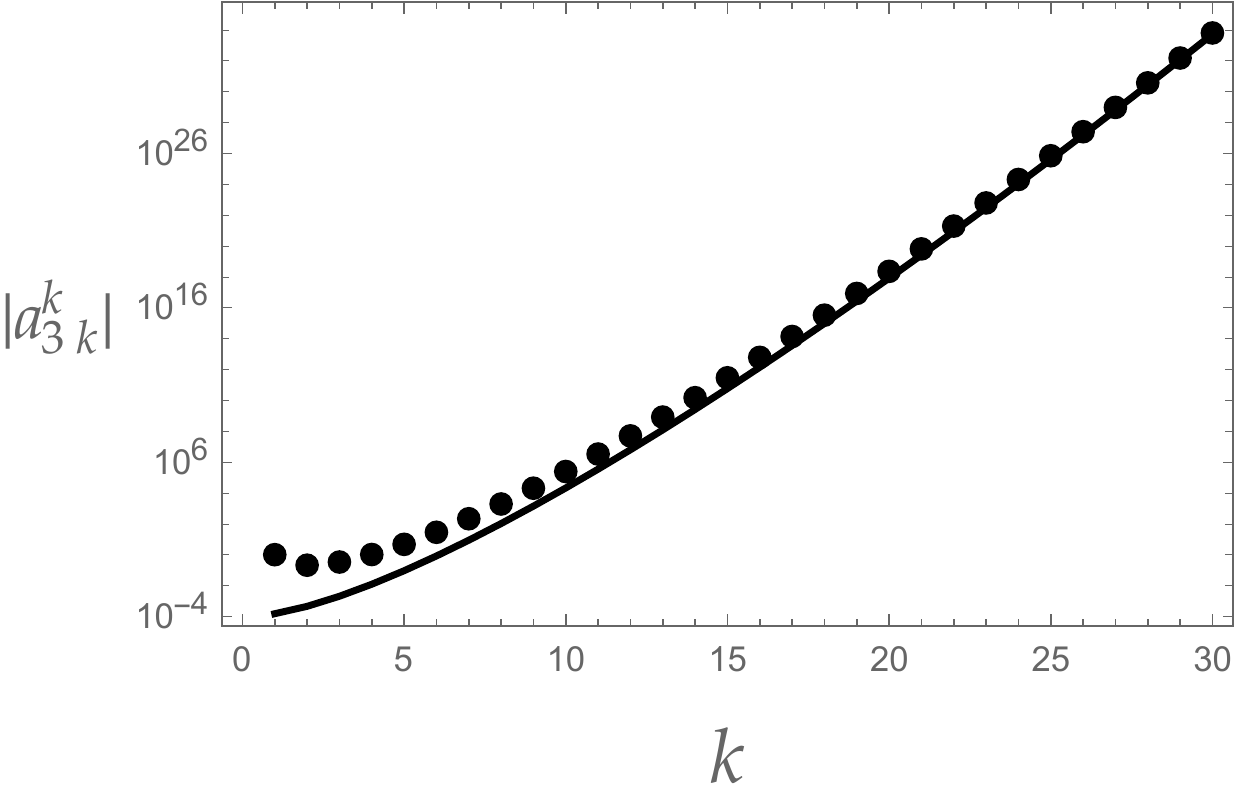}
\end{minipage}
\caption{Behaviour of $|a^{k}_{3k}|$ versus $k$ (dots), together with 
 the asymptotic law~(\ref{Eq:Debye.1.1.1.0}) evaluated with $C\,=\, 10^4$ (solid curve).}
\label{Fig:Debye.2}
\end{figure}

Such a divergence should entirely be ascribed to the asymptotic behaviour of the Debye polynomials for large values of their 
orders (for  $\epsilon < 1$ their argument is  always greater that one). Unfortunately, the mathematical properties of Debye polynomials 
are still largely unknown. On the other hand, their exploration is certainly  beyond the scopes of the present paper. 
For such reason we shall limit ourselves to grasp  the  asymptotics of $U_k(t)$ for $k\to \infty$ and  $t\gg 1$, in such a way that  
$U_k(t)\,\sim\,a^{k}_{3k}\,t^{3k}$. In particular, in Appendix~\ref{App:C} it is shown that
\begin{equation}
\label{Eq:Debye.1.1.1.0}
\begin{array}{l}
\displaystyle
a^{k}_{3k}\,\sim\,C \left(-\dfrac 32\right)^k\,k!\,,\qquad\qquad k\to \infty\,,
\end{array}
\end{equation}
where $C$ denotes a positive constant that can  be estimated empirically.
To this aim, in Fig.~\ref{Fig:Debye.2} the behaviour of $|a^{k}_{3k}|$ versus $k$ is shown (dots).
The solid curve represents the asymptotic law~(\ref{Eq:Debye.1.1.1.0}), 
evaluated with $C\,=\, 10^4$. 

\begin{figure}[!ht]
\centering
\begin{minipage}[t]{8cm}
\includegraphics[width=8cm]{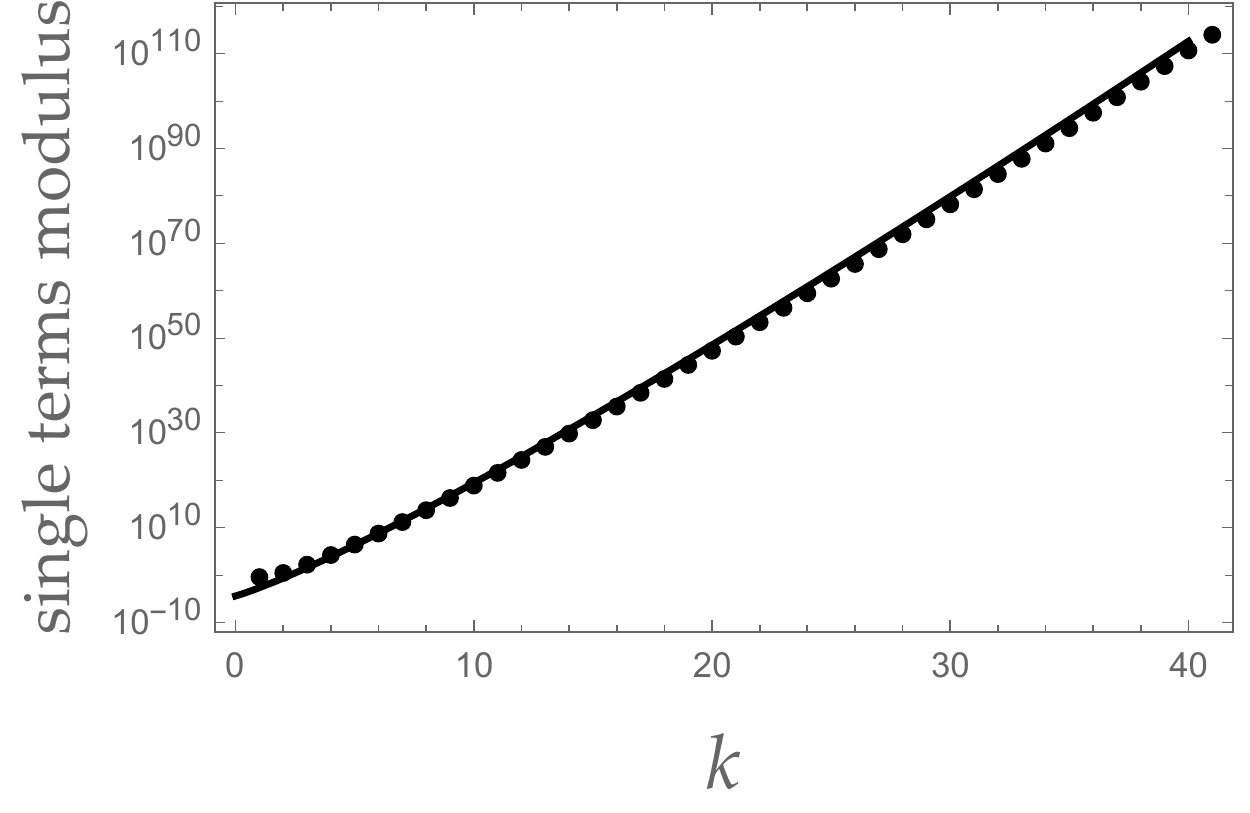}
\end{minipage}
\caption{Behaviour of the modulus of the single terms of the series~(\ref{Eq:Debye.1})
as a function of the index $k$ (dots) together with Eq.~(\ref{Eq:Debye.1.1.1.0.1}) (solid curve),
when $n=10$ and $\epsilon=99/100$.}
\label{Fig:Debye.3}
\end{figure}

The reasonably good agreement, at least at a visual level, suggests that also the single terms of the series~(\ref{Eq:Debye.1})  
could asymptotically follow, at least in the limiting case $\epsilon\to 1$, the factorial law
\begin{equation}
\label{Eq:Debye.1.1.1.0.1}
\begin{array}{l}
\displaystyle
\dfrac {\rho^n\,C}{\sqrt{2\pi\,n\,\sqrt{1\,-\,\epsilon^2}}}\,
\left(-\dfrac 3{2n (1-\epsilon^2)^{3/2}}\right)^k\,k!\,,\qquad k\to \infty\,.
\end{array}
\end{equation}
In Fig.~\ref{Fig:Debye.3} the modulus of the single terms of the series~(\ref{Eq:Debye.1})
is shown as a function of the index $k$ (dots) together with Eq.~(\ref{Eq:Debye.1.1.1.0.1}) (solid curve),
when $n=10$ and  $\epsilon=99/100$. 

Once the factorially divergent character of the expansion~(\ref{Eq:Debye.1}) has been 
established (at least on an empirical basis), it remains to see whether and how the quantity $J_n(n\epsilon)$ could be retrieved 
starting from the corresponding partial sum sequence. 
To this end, it must be recalled that several classes of factorially divergent series can be efficiently decoded by using 
suitable transformation schemes, the most known of them being 
Borel summation~\cite{Costin/2009} and Pad\'{e} approximants~\cite{Baker/Graves-Morris/1996}.
However, new resummation algorithms more effective than Borel and Pad\'e have been proposed in the last three decades for 
decoding factorially divergent alternating series. 
In the second part of our work these relatively new computational algorithms, 
called \emph{Levin-type nonlinear sequence transformations}~\cite{Levin/1973},  will  directly be applied  to 
the partial sum sequences of the series in Eq.~(\ref{Eq:Debye.1}).

For reader's convenience, in the next section a brief survey of Levin-type transformations will now be given.

\subsection{Levin-type nonlinear sequence transformations}
\label{WT}

Several nonlinear transformations  have been conceived to convert the sequence of the partial sums of a 
diverging series into new sequences converging to the so-called  \emph{generalized limit} of the series itself, namely the (finite) value ``hidden'' within the 
divergent partial sum sequence. The probably most known transformation scheme is the so-called Wynn or $\epsilon$-algorithm~\cite{Wynn/1956}, 
which implements Pad\'e approximants through a surprisingly easy recursive scheme.

Consider  the partial sum sequence, say $\{s_{n}\}_{n=0}^\infty$, of a given series $\displaystyle\sum^\infty_{k=0}\,a_k$, i.e.,
\begin{equation}
\label{Eq:Delta.1}
\begin{array}{l}
\displaystyle
s_n\,=\,\sum_{k=0}^n a_k\,,\qquad\qquad n=0,1,2,\ldots\,
\end{array}
\end{equation}
%
%
A typical feature of Wynn's algorithm (and  of several  other sequence transformations) is that only the input of the numerical values of a finite substring 
of the sequence $\{s_n\}^\infty_{n=0}$ is required to achieve resummation. In other words, no further information is needed to compute approximations 
of the limit (or the antilimit)  of~(\ref{Eq:Delta.1}), i.e., of the series sum.

However, in some cases additional information on the index dependence of the truncation errors is available. For instance, the magnitudes of the 
truncation errors of Stieltjes series are   
bounded  by the first term neglected in the partial sum  and they also display the same sign patterns. 
In 1973, Levin~\cite{Levin/1973} introduced a sequence transformation which utilizes such a valuable \emph{a priori} structural information to enhance 
the efficiency and the convergence speed of the transformed sequence.
For a general discussion about the construction principles of Levin-type sequence transformations, readers are encouraged to go through the paper by 
Weniger~\cite{Weniger/1989}  as well as through Ref.~\cite{Borghi/Weniger/2015}, which contains a slightly upgraded and formally improved version 
of the former.
%
In the same paper it was also rigorously proved that Levin-type transformations are  able to sum what is  considered to be
the paradigmatic factorially divergent series in physics, the Euler series, with convergence rates considerably higher than Pad\'e approximants.

In the following, two special Levin-type transformations will be employed to numerically explore the Stieltjes nature of the $\mathcal{S}(\epsilon; M)$.
The first is the so-called Levin $d$-transformation, which maps the partial sum sequence~(\ref{Eq:Delta.1}) into the new sequence $\{d_k\}$
defined as follows~\cite{Weniger/1989}:
\begin{equation}
\label{LevinTransformation}
    d_k=\displaystyle\frac
    {\displaystyle\sum_{j=0}^k
    \,(-1)^j \left(k\atop j\right)\,(1+j)^{k-1} 
    \displaystyle\,\frac{s_{j}}{a_{j+1}}}
    {\displaystyle
    \sum_{j=0}^k\,(-1)^j \left(k\atop j\right)\,(1+j)^{k-1} 
    \displaystyle\frac{1}{a_{j+1}}}\,,\quad\,k\,=\,1,\,2,\,\ldots
\end{equation}

The other transformation is the so-called Weniger transformation, also called $\delta$-transformation, which maps the sequence~(\ref{Eq:Delta.1}) into the sequence
$\{\delta_{k}\}_{k=1}^\infty$ defined as follows~\cite{Weniger/1989}:
\begin{equation}
\label{WenigerTransformation}
    \delta_k=\displaystyle\frac
    {\displaystyle\sum_{j=0}^k
    \,(-1)^j \left(k\atop j\right)\,(1+j)_{k-1} 
    \displaystyle\,\frac{s_{j}}{a_{j+1}}}
    {\displaystyle
    \sum_{j=0}^k\,(-1)^j \left(k\atop j\right)\,(1+j)_{k-1} 
    \displaystyle\frac{1}{a_{j+1}}}\,,\quad\,k\,=\,1,\,2,\,\ldots\,,
\end{equation}
where  $(\cdot)_{k}$ denotes Pochhammer symbol. 
It is worth nothing that in order to pass from the $\{d_k\}$ sequence to the  $\{\delta_k\}$ sequence it is sufficient
to formally replace the exponential term $(1+j)^{k-1}$ by the Pochhammer term $(1+j)_{k-1}$, 
and \emph{viceversa}. 

Despite their unusual and a first sight apparently complicated aspects, sequence transformations into 
Eqs.~(\ref{LevinTransformation}) and~(\ref{WenigerTransformation}) are of straightforward 
implementation. For instance, in the last twenty years they have profitably been employed to solve 
several optical problems, showing an extraordinary effectiveness in summing divergent as well as slowly convergent series~%
\cite{Borghi/Santarsiero/2003,borghiOL-07,borghiJOSAA-08a,borghiJOSAA-08b,borghiJOSAA-08c,borghiJOSAA-09,borghiJOSAA-9b,borghiJOSAA-10,borghiJOSAA-11,%
Borghi/Gori/Guattari/Santarsiero/2011,borghiOL-11b,borghiJOSAA-12}. 

Our first numerical experiment is aimed at appreciating the powerfulness of $d$- and $\delta$-transformations in resumming the divergent Debye expansion~(\ref{Eq:Debye.1}) 
to retrieve, from the partial sum sequences  represented in  Figs.~\ref{Fig:Debye.1}a and~\ref{Fig:Debye.1}b, the ``coded''  values of $J_{10}(5)$ and $J_{10}(9)$, respectively. 
%
Table~\ref{Tab:ConvergenceA} refers to $J_{10}(5)$: the values of the partial sum sequence $\{s_n\}$ 
of the series~(\ref{Eq:Debye.1})  are listed in the second column,
whereas the third and the fourth columns contain the values of the sequences $\{d_k\}$ and $\{\delta_k\}$  obtained by applying, on the sequence $\{s_n\}$, 
maps~(\ref{LevinTransformation}) and~(\ref{WenigerTransformation}), respectively (the corresponding orders are reported in the first column).
It can be appreciated that both sequences converge  (approximately with the same speed) toward the correct limit $J_{10}(5)=0.001467802647\ldots$
which is only partially touched by the (diverging) partial sum sequence.
\begin{table}[!ht]
\centerline{
    \begin{tabular}{|c|c|c|c|}
    \hline
  order & Partial sums 			& $d_k$ 	& $\delta_k$ 		\\ \hline \hline
  1        & \textbf{0.0014}92003408 	& -- 		& --			 	\\ \hline
 2  & \textbf{0.00146}5682591 & \textbf{0.001467}977164 & \textbf{0.001467}789214 \\ \hline
 3  & \textbf{0.00146}8263281 & \textbf{0.00146780}3250 & \textbf{0.00146780}4355 \\ \hline
 4  & \textbf{0.001467}656862 & \textbf{0.00146780}4086 & \textbf{0.001467802}513 \\ \hline
 5  & \textbf{0.0014678}63379 & \textbf{0.001467802}576 & \textbf{0.0014678026}31 \\ \hline
 6  & \textbf{0.001467}770986 & \textbf{0.0014678026}34 & \textbf{0.00146780264}1 \\ \hline
 7  & \textbf{0.0014678}22501 & \textbf{0.00146780264}2 & \textbf{0.00146780264}6 \\ \hline
 8  & \textbf{0.001467}788088 & \textbf{0.001467802647} & \textbf{0.001467802647} \\ \hline
 9  & \textbf{0.001467}814880 & \textbf{0.001467802647} & \textbf{0.001467802647} \\ \hline
 10  & \textbf{0.001467}791058 & \textbf{0.001467802647} & \textbf{0.001467802647} \\ \hline
 11  & \textbf{0.0014678}14875 & \textbf{0.001467802647} & \textbf{0.001467802647} \\ \hline
 12  & \textbf{0.001467}788427 & \textbf{0.001467802647} & \textbf{0.001467802647} \\ \hline
 13  & \textbf{0.001467}820725 & \textbf{0.001467802647} & \textbf{0.001467802647} \\ \hline
 14  & \textbf{0.001467}777704 & \textbf{0.001467802647} & \textbf{0.001467802647} \\ \hline
 15  & \textbf{0.0014678}39774 & \textbf{0.001467802647} & \textbf{0.001467802647} \\ \hline
 16  & \textbf{0.001467}743345 & \textbf{0.001467802647} & \textbf{0.001467802647} \\ \hline
 17  & \textbf{0.001467}903836 & \textbf{0.001467802647} & \textbf{0.001467802647} \\ \hline
 18  & \textbf{0.001467}618939 & \textbf{0.001467802647} & \textbf{0.001467802647} \\ \hline
 19  & \textbf{0.00146}8156252 & \textbf{0.001467802647} & \textbf{0.001467802647} \\ \hline
 20  & \textbf{0.001467}083337 & \textbf{0.001467802647} & \textbf{0.001467802647} \\ \hline
 21  & \textbf{0.00146}9344663 & \textbf{0.001467802647} & \textbf{0.001467802647} \\ \hline
 22  & \textbf{0.00146}4327946 & \textbf{0.001467802647} & \textbf{0.001467802647} \\ \hline
 23  & \textbf{0.0014}76013517 & \textbf{0.001467802647} & \textbf{0.001467802647} \\ \hline
 24  & \textbf{0.0014}47498723 & \textbf{0.001467802647} & \textbf{0.001467802647} \\ \hline
 25  & \textbf{0.001}520240513 & \textbf{0.001467802647} & \textbf{0.001467802647} \\ \hline
 26  & \textbf{0.001}326611312 & \textbf{0.001467802647} & \textbf{0.001467802647} \\ \hline
 27  & \textbf{0.001}863491524 & \textbf{0.001467802647} & \textbf{0.001467802647} \\ \hline
 28  & 0.0003153551669 & \textbf{0.001467802647} & \textbf{0.001467802647} \\ \hline
 29  & 0.004951150350 & \textbf{0.001467802647} & \textbf{0.001467802647} \\ \hline
 30  & -0.009444360750 & \textbf{0.001467802647} & \textbf{0.001467802647} \\ \hline
  \ldots & \ldots  & \ldots & \ldots  \\ \hline
    \hline
    \end{tabular}
}
\caption{Using Levin-type transformations to retrieve the value of $J_{10}(5)$ via Debye's expansion~(\ref{Eq:Debye.1}).
2nd column: partial sum sequence. 3rd column: $d$-transformation sequence. 4th column: $\delta$-transformation sequence.
The correct limit is $J_{10}(5)=0.001467802647\ldots$}
\label{Tab:ConvergenceA}
\end{table}

Things are much  more evident as far as the evaluation of $J_{10}(9)$ is concerned. 
This is shown in Table~\ref{Tab:ConvergenceB}, where it is easier to grasp the alternating divergent 
character of the partial sum sequence, dramatically evident since order 5. 
At the same time, it is now possible to appreciate the different convergence speed of $d$- and $\delta$-transformations
in retrieving the correct value $J_{10}(9)=0.1246940928\ldots$
\begin{table}[!ht]
\centerline{
    \begin{tabular}{|c|c|c|c|}
    \hline
  order & Partial sums  				& $d_k$ 	& $\delta_k$ 		\\ \hline \hline
 1  & \textbf{0.1}397916170 		& -- & --  \\ \hline
 2  & \textbf{0.1}086355082 		& \textbf{0.12}54181699 & \textbf{0.124}0036791  \\ \hline
 3  & \textbf{0.1}617358916		& \textbf{0.124}8610123 & \textbf{0.124}6183759  \\ \hline
 4  & -0.01216322740 			& \textbf{0.1246}749430 & \textbf{0.124}7070912  \\ \hline
 5  & 0.8321487102 				& \textbf{0.12469}95597 & \textbf{0.124}7020129  \\ \hline
 6  & -4.608269328 				& \textbf{0.124}7001707 & \textbf{0.12469}65877  \\ \hline
 7  & 39.11010231 				& \textbf{0.12469}52850 & \textbf{0.124694}8554  \\ \hline
 8  & -381.9081096 				& \textbf{0.124694}3011 & \textbf{0.124694}3936  \\ \hline
 9  & 4344.426282 				& \textbf{0.124694}2503 & \textbf{0.124694}2448  \\ \hline
 10  & -56259.36907 				&\textbf{ 0.124694}1463 & \textbf{0.124694}1756  \\ \hline
 11  & 817636.3501 				& \textbf{0.12469409}39 & \textbf{0.124694}1370  \\ \hline
 12  & -1.317999719 $\times 10^7$ 	& \textbf{0.1246940}899 & \textbf{0.124694}1153  \\ \hline
 13  & 2.333958899  $\times 10^8$ 	& \textbf{0.124694092}0 & \textbf{0.124694}1037  \\ \hline
 14  & -4.504271888 $\times 10^9$ 	& \textbf{0.124694092}1 & \textbf{0.12469409}78  \\ \hline
 15  & 9.409762678 $\times 10^{10}$ & \textbf{0.124694092}3 & \textbf{0.12469409}48  \\ \hline
 16  & -2.11566839 3$\times 10^{12}$ & \textbf{0.124694092}6 & \textbf{0.12469409}35  \\ \hline
 17  & 5.094033071 $\times 10^{13}$  & \textbf{0.1246940928} & \textbf{0.124694092}9  \\ \hline
 18  & -1.307753975 $\times 10^{15}$ & \textbf{0.1246940928} & \textbf{0.124694092}6  \\ \hline
 19  & 3.565916241 $\times 10^{16} $ & \textbf{0.1246940928} & \textbf{0.124694092}6  \\ \hline
 20  & -1.029237477 $\times 10^{18}$ & \textbf{0.1246940928} & \textbf{0.124694092}6  \\ \hline
 21  & 3.134988579 $\times 10^{19}$ & \textbf{0.1246940928} & \textbf{0.124694092}6  \\ \hline
 22  & -1.004946391 $\times 10^{21}$ & \textbf{0.1246940928} & \textbf{0.124694092}7  \\ \hline
 23  & 3.381908041 $\times 10^{22}$ & \textbf{0.1246940928} & \textbf{0.124694092}7  \\ \hline
 24  & -1.192111934 $\times 10^{24}$ & \textbf{0.1246940928} & \textbf{0.124694092}7  \\ \hline
 25  & 4.392572423 $\times 10^{25}$ & \textbf{0.1246940928} & \textbf{0.1246940928}  \\ \hline
   \ldots & \ldots  & \ldots & \ldots  \\ \hline
   \hline
    \end{tabular}
}
\caption{The same as in Tab.~\ref{Tab:ConvergenceA} but for  $J_{10}(9)=0.1246940928\ldots$
}
\label{Tab:ConvergenceB}
\end{table}

\subsection{Is Kepler Kapteyn series Stieltjes?}
\label{Sec:KKStieltjes}

Equipped with the above mathematical tools (Debye expansion and Levin-type transformations),
we are now ready to explore the Stieltjes nature of the KS~(\ref{Eq:KE.4}).
To this end,  Eq.~(\ref{Eq:Debye.1}) is first substituted into Eq.~(\ref{Eq:KE.4}) to give
\begin{equation}
\label{Eq:Debye.1.1.1}
\begin{array}{l}
\displaystyle
\S(\epsilon;M)\,=\,\sqrt{\dfrac{2}{\pi\sqrt{1-\epsilon^2}}}\,\sum_{n=1}^\infty\,\sum_{k=0}^\infty\,
\dfrac {z^n}{n^{k+3/2}}\,U_k\,\left(\dfrac 1{\sqrt{1-\epsilon^2}}\right)\,,
\end{array}
\end{equation}
which, on interchanging the order of the $n$- and of the $k$- series, leads to
\begin{equation}
\label{Eq:Debye.1.1.2}
\begin{array}{l}
\displaystyle
\S(\epsilon;M)\,=\,
\sqrt{\dfrac{2}{\pi\sqrt{1-\epsilon^2}}}\,\sum_{k=0}^\infty\,\,U_k\,\left(\dfrac 1{\sqrt{1-\epsilon^2}}\right)\,
\sum_{n=1}^\infty\,\dfrac {z^n}{n^{k+3/2}}\,,
\end{array}
\end{equation}
and finally, on taking Eq.~(\ref{Eq:KE.6.1}) into account, to
\begin{equation}
\label{Eq:Debye.1.1.2.1}
\begin{array}{l}
\displaystyle
\S(\epsilon;M)\,=\,
\displaystyle
\sqrt{\dfrac{2}{\pi\sqrt{1-\epsilon^2}}}\,
\sum_{k=0}^\infty\,U_k\,\left(\dfrac 1{\sqrt{1-\epsilon^2}}\right)\,\L_{k+3/2}(z)\,.
\end{array}
\end{equation}
The subsequent step consists in substituting the integral representation of polylogarithm function given by  
Eq.~(\ref{Eq:Debye.1.1.2.3}) into Eq.~(\ref{Eq:Debye.1.1.2.1}) and then to interchange the series with the integral. 
After simple algebra the following integral representation of 
$\S(\epsilon,M)$ can then be established:
\begin{equation}
\label{Eq:Debye.1.1.2.5}
\begin{array}{l}
\displaystyle
\S(\epsilon;M)\,=\,
\displaystyle
-\sqrt{\dfrac{2}{\pi\sqrt{1-\epsilon^2}}}\,
\int_0^1\,
\dfrac{\mathcal{U}\left(-\log t,\dfrac 1{\sqrt{1-\epsilon^2}}\right)}{t\,-\,1/z}\,\d t\,,
\end{array}
\end{equation}
with the function $\mathcal{U}(\cdot,\cdot)$ being defined as follows:
\begin{equation}
\label{Eq:Debye.1.1.2.6}
\begin{array}{l}
\displaystyle
\mathcal{U}\left(x,y\right)\,=\,
\sum_{k=0}^\infty\,
\dfrac{x^{k+1/2}}{\Gamma\left(k+\dfrac 32\right)}\,U_k(y)\,,\quad x\ge 0,\,y > 1.
\end{array}
\end{equation}
%
%
In order to prove the ``Stieltjness'' of  $\mathcal{S}$, it would be mandatory to show that, 
for any given $y > 1$, both $\mathcal{U}$ and  $\partial \mathcal{U}/\partial x$ to be nonnegative for $x\in [0,\infty)$.
Presently we do not possess a rigorous proof of it, so that we shall limit ourselves to provide strong numerical evidences of such a conjecture.

In doing so, Levin-type transformations will be again employed to decode the series into Eq.~(\ref{Eq:Debye.1.1.2.6}) which,
thanks to the above described  asymptotics of Debye's polynomials, is expected to display a factorial divergence.
The result of a preliminary numerical experiment is shown in Tab.~\ref{Tab:ConvergenceC}, where $d$ and $\delta$  were employed to 
retrieve the value of $\mathcal{U} \left(\log 2,\dfrac {100}{\sqrt{199}}\right)$,  
corresponding to the pair $(t,\epsilon)=(1/2,{99}/{100})$.  
\begin{table}[!ht]
\centerline{
    \begin{tabular}{|c|c|c|c|}
    \hline
  order & Partial sums  			& $d_k$ 	& $\delta_k$ 		\\ \hline \hline
1  & 0.9394372787 & -- & -- \\
 2  & -30.89260169 & 0.08352018113 & 0.5787695135 \\
 3  & 4951.945127 & 0.1804392091 & 0.5112075442 \\
 5  & 2.620274608  $\times 10^8$ & 0.3855265022 & \textbf{0.4}507501433 \\
 10  & -5.444869076  $\times 10^{20}$ & \textbf{0.41}23795594 & \textbf{0.41}45089856 \\
 15  & 1.878304379  $\times 10^{33}$ & \textbf{0.41}31217162 & \textbf{0.41}19027710 \\
 20  & -7.870085134  $\times 10^{45}$ & \textbf{0.412874}4275 & \textbf{0.412}5618326 \\
 25  & 3.658039660  $\times 10^{58}$ & \textbf{0.4128567}548 & \textbf{0.412}9445672 \\
 30  & -1.813835802  $\times 10^{71}$ & \textbf{0.412857}3130 & \textbf{0.41}30283827 \\
 35  & 9.400315017  $\times 10^{83}$ & \textbf{0.41285746}19 & \textbf{0.41}30041720 \\
 40  & -5.030871012  $\times 10^{96}$ & \textbf{0.41285746}59 & \textbf{0.412}9603691 \\
 45  & 2.758995841  $\times 10^{109}$ & \textbf{0.412857464}9 & \textbf{0.412}9237187 \\
 50  & -1.542390197  $\times 10^{122}$ & \textbf{0.4128574648} & \textbf{0.4128}982813 \\
 55  & 8.757110192  $\times 10^{134}$ & \textbf{0.4128574648} & \textbf{0.4128}819512 \\
 60  & -5.035764309  $\times 10^{147}$ & \textbf{0.4128574648} & \textbf{0.4128}718730 \\
 65  & 2.926920321  $\times 10^{160}$ & \textbf{0.4128574648} & \textbf{0.4128}657968 \\
 70  & -1.716729904  $\times 10^{173}$ & \textbf{0.4128574648} & \textbf{0.4128}621941 \\
 75  & 1.014819268  $\times 10^{186}$ & \textbf{0.4128574648} & \textbf{0.4128}600897 \\
 80  & -6.039885436  $\times 10^{198}$ & \textbf{0.4128574648} & \textbf{0.41285}88796 \\
 85  & 3.616269264  $\times 10^{211}$ & \textbf{0.4128574648} & \textbf{0.41285}81966 \\
 90  & -2.176637686  $\times 10^{224}$ & \textbf{0.4128574648} & \textbf{0.412857}8202 \\
 95  & 1.316300235  $\times 10^{237}$ & \textbf{0.4128574648} & \textbf{0.412857}6192 \\
 100  & -7.993851066  $\times 10^{249}$ & \textbf{0.4128574648} & \textbf{0.412857}5168 \\
 105  & 4.873145754  $\times 10^{262}$ & \textbf{0.4128574648} & \textbf{0.41285746}83 \\
   \ldots & \ldots  & \ldots & \ldots  \\ \hline
   \hline
    \end{tabular}
}
\caption{Using Levin-type transformations to retrieve the value of $\mathcal{U} \left(\log 2,\dfrac {100}{\sqrt{199}}\right)$ via 
the divergent expansion given in Eq.~(\ref{Eq:Debye.1.1.2.6}).
}
\label{Tab:ConvergenceC}
\end{table}
Also in that case the convergence of both $d$-transformed and $\delta$-transformed sequences is more than evident.
However, the performances of the $d$-transformation appears to be slightly superior with respect those of $\delta$. Such a 
superiority has been confirmed during other several tests carried out on different pairs $(t,\epsilon)$ 
(not shown here for evident room reasons). For such reason, in the rest of the paper only  the 
Levin $d$-transformation will be used to provide evidences of the ``Stieltjness conjecture''. 
%
\begin{figure}[!ht]
\centering
\begin{minipage}[t]{8cm}
\includegraphics[width=8cm]{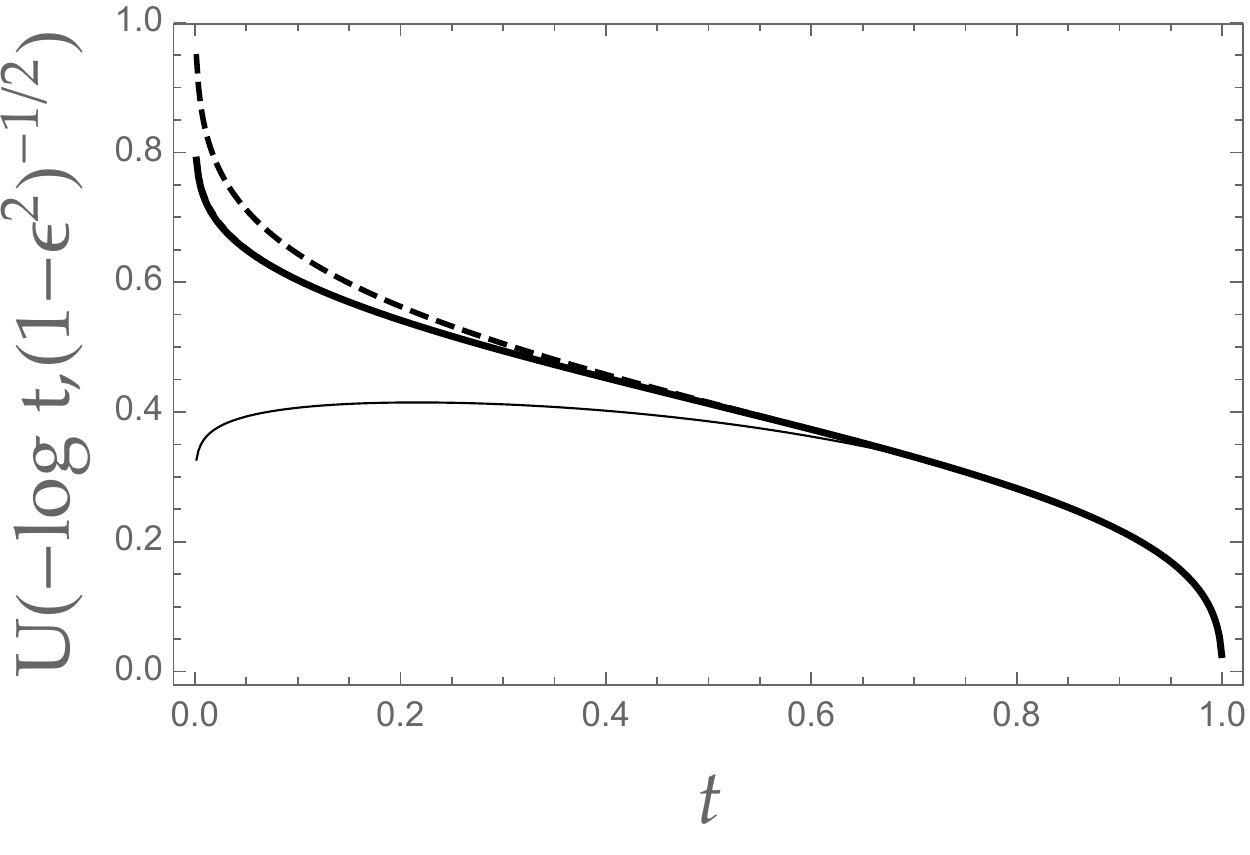}
\end{minipage}
\caption{Behaviour of $\mathcal{U}(-\log t, 1/{\sqrt{1-\epsilon^2}})$ against $t$ for different values of the 
$d$-transformation order, precisely 6~(solid thin curve), 10~(dashed curve), 20~(dotted curve), 
and  40~(thick solid curve). The orbit eccentricity is $\epsilon=99/100$.}
\label{Fig:Debye.3.99}
\end{figure}

To this end, in Fig.~\ref{Fig:Debye.3.99} the behaviour of $\mathcal{U}\left(-\log t,\dfrac 1{\sqrt{1-\epsilon^2}}\right)$ has been plotted for $t\in (0,1]$ and for $\epsilon=99/100$. 
The four curves correspond to  different values of the $d$-transformation order,  precisely 
6~(solid thin curve), 10~(dashed curve), 20~(dotted curve), 
and  40~(thick solid curve). The last two curve are practically coincident, which confirms the rapid
convergence of $d$-transformation for all  $t \in (0,1]$.
More importantly, the {\em shape} of the graphic of the generating  function $\mathcal{U}$ seems to possess those  
features which are mandatory for integral~(\ref{Eq:Debye.1.1.2.5}) to be a \emph{Stieltjes integral}.
\begin{figure}[!ht]
\centering
\begin{minipage}[t]{8cm}
\includegraphics[width=8cm]{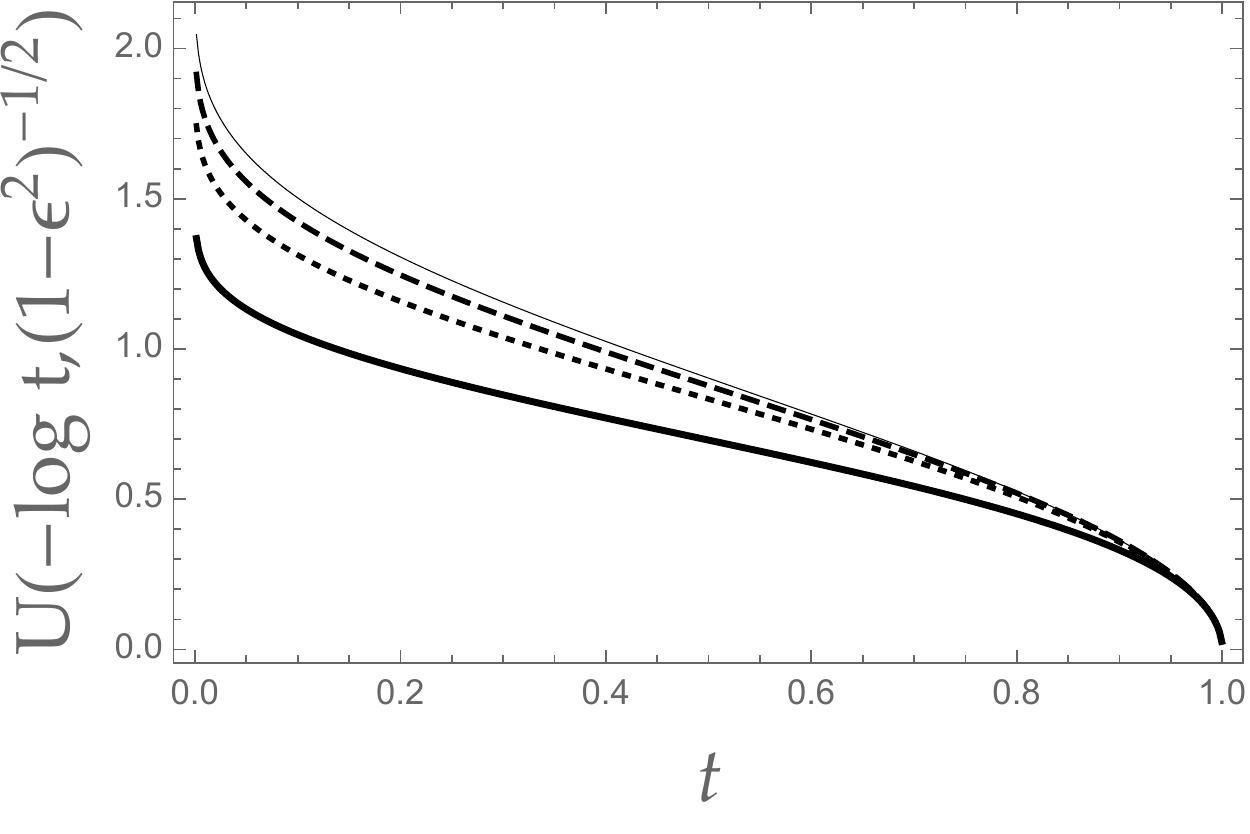}
\end{minipage}
\caption{Behaviour of $\mathcal{U}(-\log t, 1/{\sqrt{1-\epsilon^2}})$ against $t$ for 
for $\epsilon=1/10$ (solid thin curve),  $\epsilon=5/10$ (dashed curve), $\epsilon=7/10$ (dotted curve), 
and $\epsilon=9/10$ (solid thick curve).  
All curves have been generated by setting  the $d$-transformation  order at 40. 
}
\label{Fig:Debye.4}
\end{figure}

Smaller values  of $\epsilon$ are expected to give much less numerical troubles that 99/100.
In Fig.~\ref{Fig:Debye.4} the behaviour of $\mathcal{U}\left(-\log t,\dfrac 1{\sqrt{1-\epsilon^2}}\right)$ is 
still represented vs. $t$ for $\epsilon=1/10$ (solid thin curve),  $\epsilon=5/10$ (dashed curve), $
\epsilon=7/10$ (dotted curve), and $\epsilon=9/10$ (solid thick curve).  
All curves have been generated by setting  the $d$-transformation  order at 40, a value which guarantees 
a rapid convergence for any $t$, as confirmed by the results of several numerical tests not reported here.
Moreover, further numerical experiments have also been carried out for different values of 
$\epsilon$ and of the  $d$-transformation order.

Figures~\ref{Fig:Debye.3.99} and~\ref{Fig:Debye.4} are  the main results of the present paper. 
%
%
They strongly  support our conjecture about the ``Stieltjesness'' of the 
complex  KS series representation of  $\S$  in Eq.~(\ref{Eq:KE.4}). 


\section{Solving Kepler's equation via Levin-type transformations}
\label{Sec:KeplerWT}

As it was said before, nonlinear sequence transformations in general and Levin-type transformation  in 
particular have proved to have extraordinary capabilities as far as Stieltjes asymptotic series decoding is 
concerned. However,  differently from Pad\'e approximants, at present a convergence theory of 
nonlinear sequence transformations on Stieltjes series  is far from being complete. 
Accordingly, all numerical experiments within such a perspective should be welcomed.
In the previous sections the Stieltjes nature of the function $\S$  has been conjectured
and numerically ``validated'' on the basis of several experiments. 
It is now time to analyze the performances of  $d$- and $\delta$-transformations when they are directly applied 
on the KS series in Eq.~(\ref{Eq:KE.4}) to solve Kepler's Equation~(\ref{Eq:KE.1}).

\begin{figure}[!ht]
\centering
\begin{minipage}[t]{8cm}
\includegraphics[width=6cm]{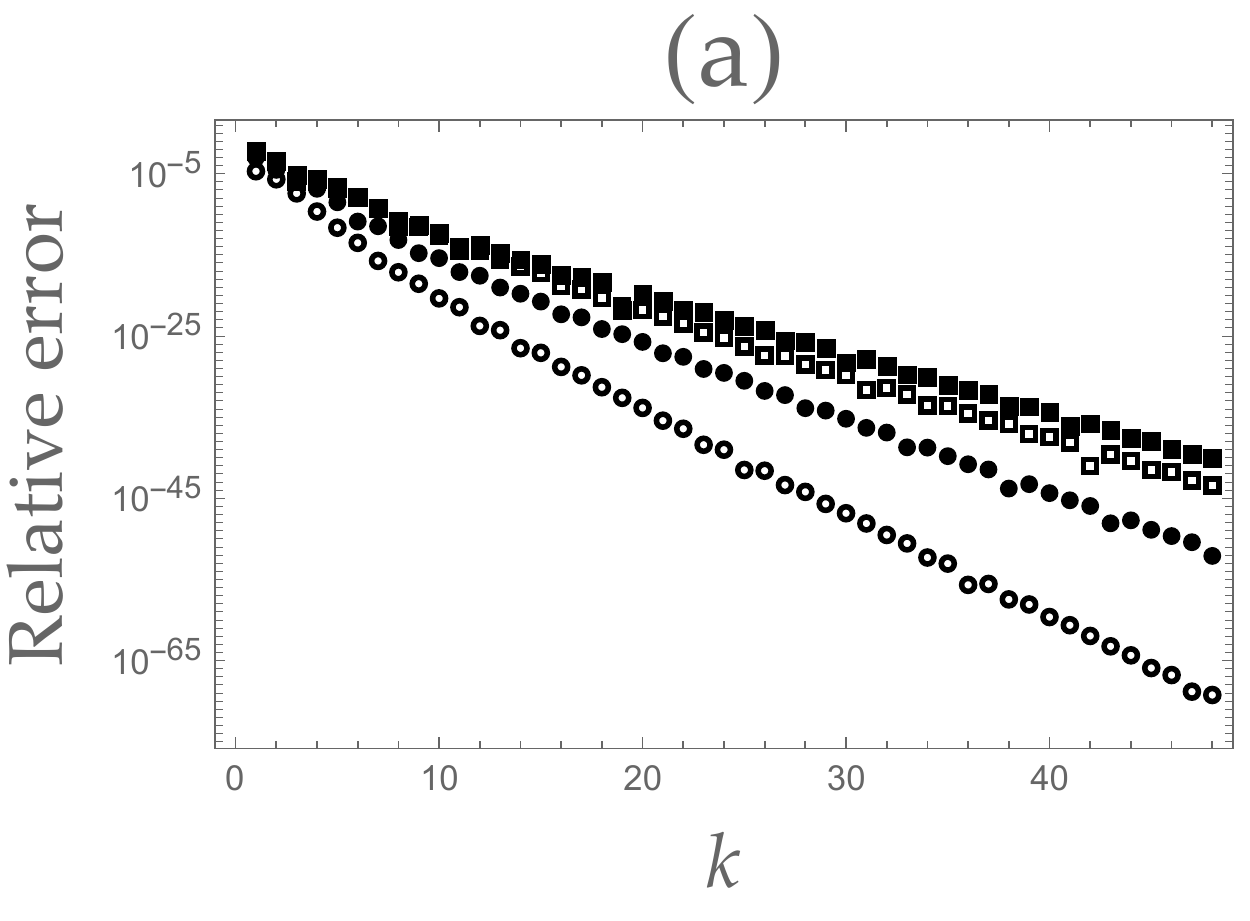}
\includegraphics[width=6cm]{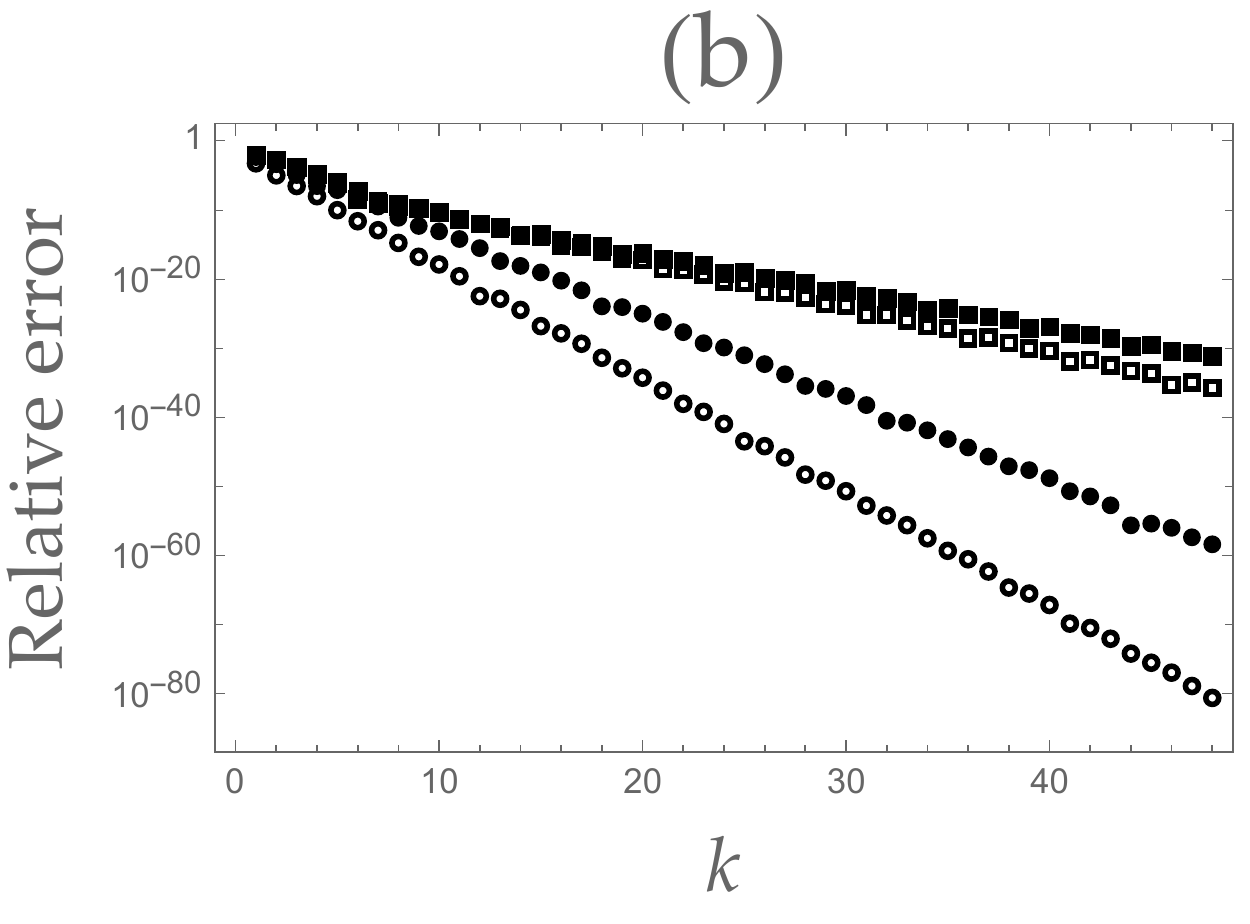}
\end{minipage}
\caption{
Plot of the relative error, evaluated through Eqs.(\ref{LevinTransformation}) and~(\ref{WenigerTransformation}) with 
respect to the value obtained by directly solving Eq.~(\ref{Eq:KE.1}) via standard numerical 
methods, against the $\delta$- (a) and the $d$- (b) transformation order $k$. 
$M=\pi/2$,  $\epsilon=2/10$ (open circles), $\epsilon=6/10$ 
(black circles), $\epsilon=9/10$ (open squares), and $\epsilon=99/100$ (black squares).}
\label{Fig:KeplerWT.1}
\end{figure}

Consider as a first example the solution of KE for $M=\pi/2$ and 
for different values of $\epsilon$. In Fig.~\ref{Fig:KeplerWT.1} the relative error, evaluated 
through Eqs.(\ref{LevinTransformation}) and~(\ref{WenigerTransformation}), with respect to the value obtained by directly solving 
Eq.~(\ref{Eq:KE.1}) via standard numerical methods, against the $\delta$- (a) and the 
$d$- (b) transformation order $k$. The open circles corresponds to $\epsilon=2/10$, the black circles to 
$\epsilon=6/10$, the open squares to $\epsilon=9/10$, and the black squares to $\epsilon=99/100$.

\begin{figure}[!ht]
\centering
\begin{minipage}[t]{8cm}
\includegraphics[width=8cm]{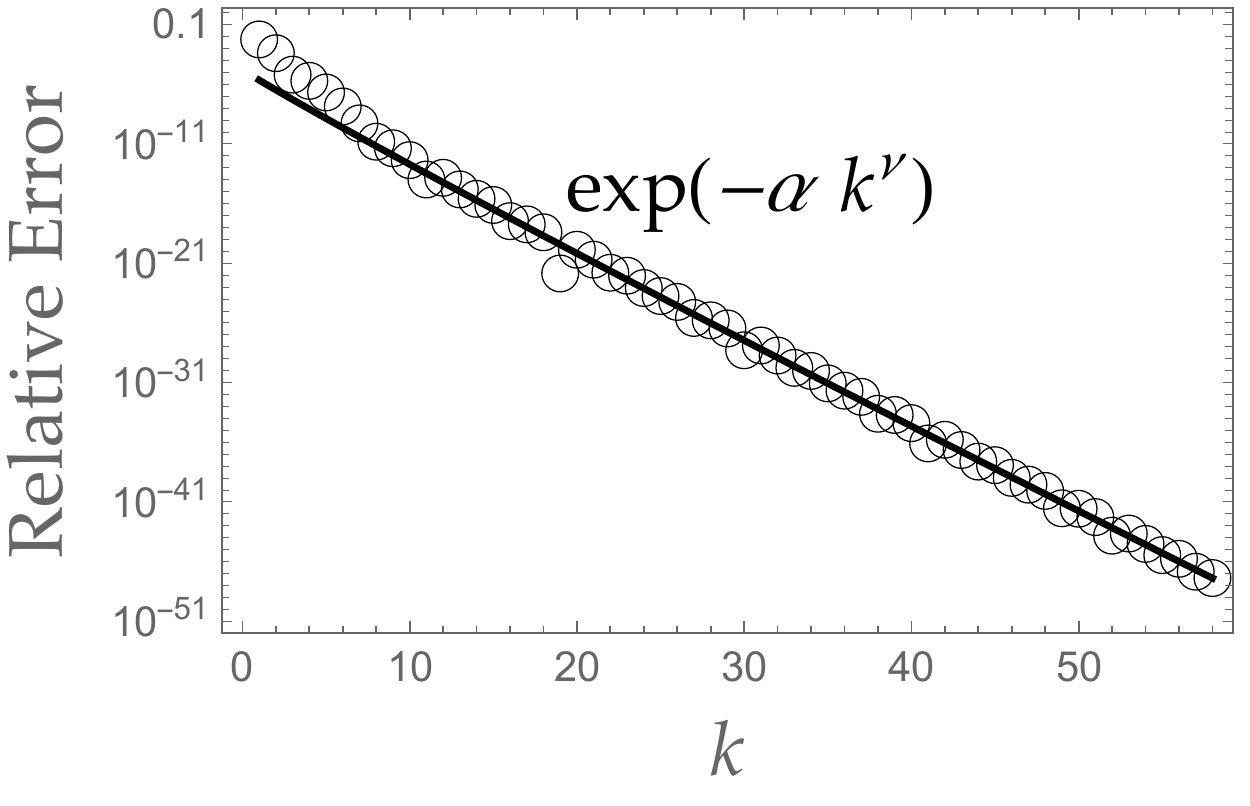}
\end{minipage}
\caption{
Relative transformation errors of Weniger's $\delta$ 
transformation~(\ref{WenigerTransformation}) vs the transformation order $k$ (open circles)
for $\epsilon=99/100$ and $M=\pi/2$.
The solid curve corresponds to a purely numerical fit of the exponential model 
$\exp(-\alpha k^\nu)$, where $\nu \simeq 1$.
}
\label{Fig:KeplerWT.1.1}
\end{figure}

In order to extract a quantitative information about the convergence rates of the $\delta$-transformation, in 
Fig.~\ref{Fig:KeplerWT.1.1} the behaviour of the relative transformation error is plotted,
as a function of the transformation order  $k$, for $M=\pi/2$ and $\epsilon=99/100$ 
(open circles). The solid curve represents the analytical fit $\exp(-\alpha k^\nu)$, where 
the parameter $\nu$ provides a meaningful measure of the convergence rate.

\begin{figure}[!ht]
\centering
\begin{minipage}[t]{8cm}
\includegraphics[width=8cm]{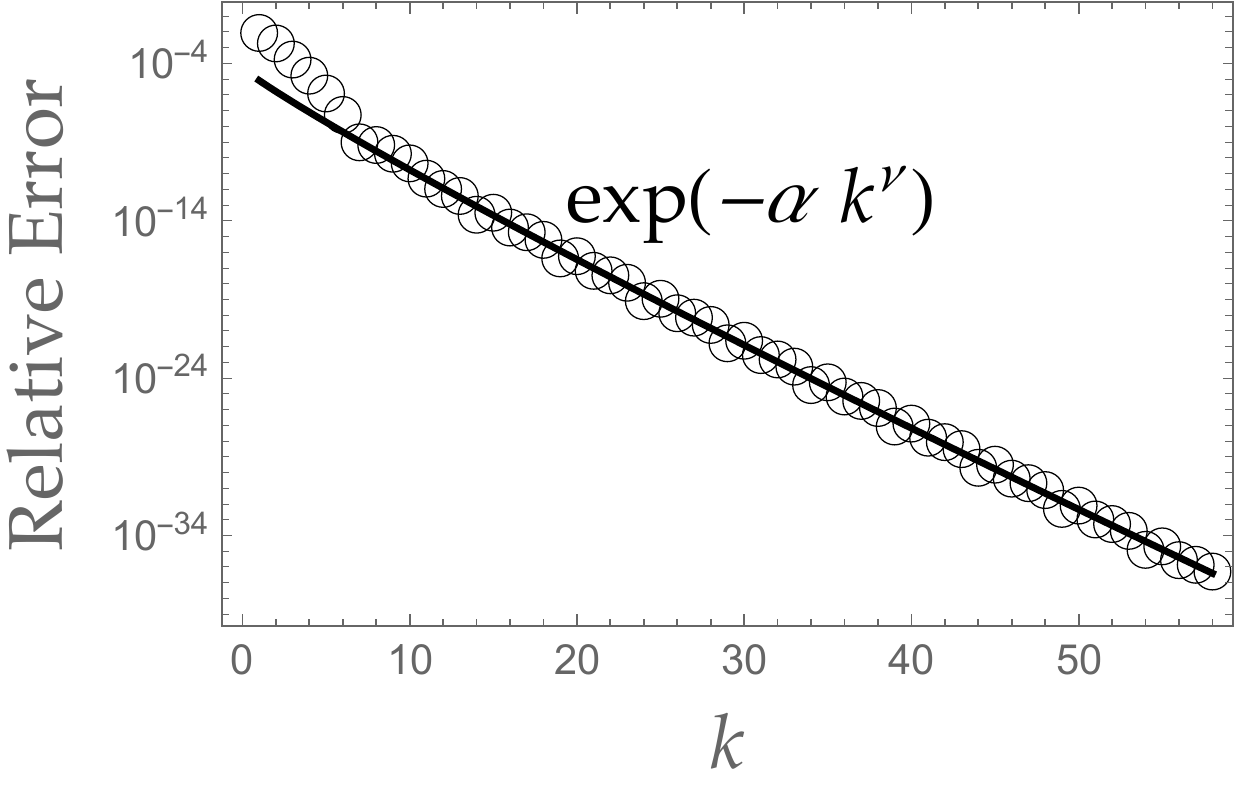}\end{minipage}
\caption{
The same as in Fig.~\ref{Fig:KeplerWT.1.1} but for the Levin's $d$ 
transformation~(\ref{LevinTransformation}). Our results indicate $\nu \simeq 9/10$.
}
\label{Fig:KeplerWT.1.1.1}
\end{figure}

A similar analysis was carried out in~\cite{Borghi/Weniger/2015} to estimate the convergence 
rate, in the limit of $k\gg 1$, of both  $\delta$- and  $d$-transformations, together with Pad\'e approximants, 
when they were applied to the Euler series. It was there found that $\nu= 2/3$ for the 
$\delta$-transformation and $\nu\simeq 3/4$ for the $d$-transformation,
both of them greater than the convergence rate of Pad\'e approximants, $\nu=1/2$~\cite{Borghi/Weniger/2015}. For 
the Kepler equation  a similar analysis can  be carried out, although necessarily from a numerical perspective, as it 
is shown for example  in Fig.~\ref{Fig:KeplerWT.1.1}, corresponding to the case $M=\pi/2$ and $\epsilon=99/100$. 
The fit is evaluated on the tail (i.e., for $k>10$), and gives $\nu\simeq 1$.
Similarly as it was done in~\cite{Borghi/Weniger/2015}, Fig.~\ref{Fig:KeplerWT.1.1.1}
shows the same as in Fig.~\ref{Fig:KeplerWT.1.1} but for the Levin $d$-transformation
in place of the $\delta$-transformation. 

\begin{figure}[!ht]
\centering
\begin{minipage}[t]{8cm}
\centerline{\includegraphics[width=8cm]{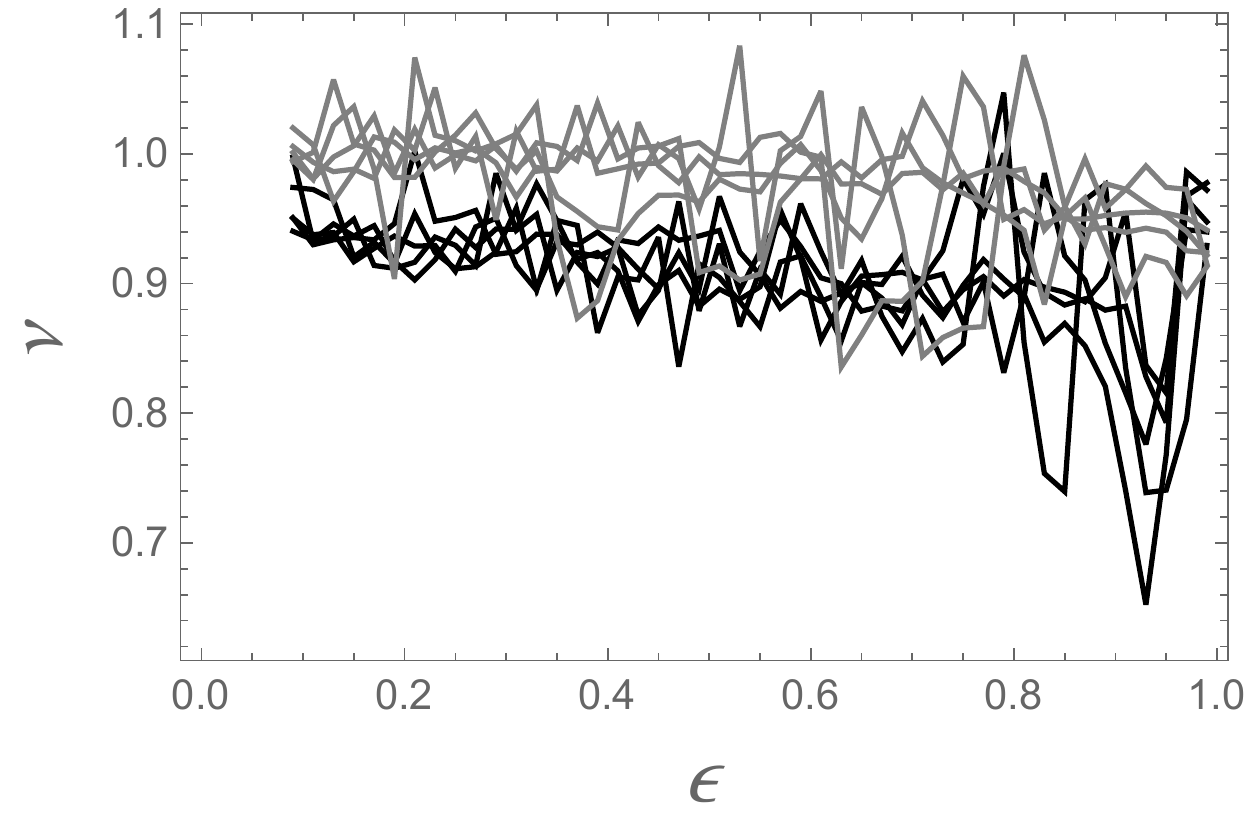}}
\end{minipage}
\caption{Behaviour of the convergence parameter $\nu$ vs. $\epsilon$ for 
$M=\pi/4\,,\,\pi/3\,,\,\pi/2\,,\,2\pi/3\,,\,3\pi/4$.
Black solid curves: $\delta$-transformation. Gray solid curves: $d$-transformation.}
\label{Fig:KeplerWT.1.2}
\end{figure}

In the present case the convergence rate of $\delta$-transformation is slightly greater than that of $d$-transformation.
However, a more detailed numerical investigation reveals that the latter is on average more
efficient than the former, as reported in Fig.~\ref{Fig:KeplerWT.1.2}, where the 
behaviour of $\nu$ is plotted against $\epsilon\in [1/10,99/100]$, for
$M=\pi/4\,,\,\pi/3\,,\,\pi/2\,,\,2\pi/3\,,\,3\pi/4$. Black  curves correspond to Weniger
$\delta$-transformation, gray curves to  Levin $d$-transformation.

\section{Discussion and conclusions}
\label{Sec:Conclusions}

The classical Bessel solution of Kepler's equation has here been revisited 
from a new perspective, offered by the still largely unexplored framework 
of the nonlinear sequence transformations.

After introducing a complex version of the original  Kapteyn series~(\ref{Eq:KE.3}) representing KE's solution, numerical evidences supporting the Stieltjes conjecture for such Kapteyn series were provided. 
Our principal tool was the (divergent) Debye expansion of Bessel functions~(\ref{Eq:Debye.1}), whose numerical 
exploration was carried out by using two different Levin-type nonlinear sequence transformations, the 
Levin $d$ and the Weniger $\delta$.
As a potentially interesting byproduct of our analysis, an effective recursive algorithm to generate arbitrarily 
high-order Debye's polynomials has also been developed. 

Exploring the Stieltjes conjecture for the Kepler equation Bessel solution could also
reveal new relevant aspects concerning the more general framework of Kapteyn series,
which continues to play a fundamental role in several areas of theoretical physics and mathematics. 
The need of developing new methods to achieve  summability of Kapteyn series of  both first and second kind
have already been invoked in the past~\cite{Tautz/Lerche/Dominici/2011}. 
In this respect, we believe the results presented here could be of some help.
In fact,  it would be rather natural ask whether Levin-type transformations would still 
be so effective in resumming more general type  of Kapteyn series, similarly as 
done in~\cite{Tautz/Lerche/Dominici/2011} on using Wynn's $\epsilon$-algorithm.
Just to give an idea, consider what follows.

Several numerical tests (not shown here) carried out on
the Kapteyn series $\displaystyle\sum_{m\ge 1}\,\dfrac{z^m}m\,J_m(m\epsilon)$,
showed Levin-type transformations to be able in resumming  the series not only for $|z|=1$, which corresponds to the 
KE case, but also for $|z| > 1$, provided that $|\mathrm{arg}(z)|< \pi$. 
This could be a direct consequence of the Stieltjes conjecture. In other words, 
the Kapteyn series  could be ``numerically continued'' to a function $F(z;\epsilon)$ living, 
for $\epsilon \in [0,1)$, on the whole complex plane but the infinite interval $(1,\infty)$.
However, things are more complex than they could appear. 
In Tab.~\ref{Tab.Finale.KS} the performances of Levin $d$- and of Weniger $\delta$-transformations in resumming the Kapteyn series for 
$\epsilon=9/10$ and $z= 10\,\exp(\mathrm{i}\pi/3)$ are shown. The indisputable superiority of $\delta$ over $d$ is now more than evident, 
the former being quickly convergent toward a precise limit. But, is such a limit meaningful?
And if so, what should be its meaning?
\begin{table}[!ht]
\centerline{
    \begin{tabular}{|c|c|c|c|}
    \hline
  order & Partial sums 			& $d_k$ 	& $\delta_k$ 		\\ \hline \hline
 1 & 2.02+3.51 i & -1.159850 + 0.307107 i & 0.112240 + 1.211289 i \\
 10 & (4.4 - 10. i) $\times 10^8$ 
 							& -1.000290 + 1.238221 i & -1.003096 +1.238166 i \\
 20 & (-3.1 + 32. i)  $\times 10^{18}$ 
 							& -1.001697 + 1.238760 i & -1.001839 + 1.238763 i \\
 30 & (7.7 +10. i) $\times 10^{27}$ 
 							& -1.001977 + 1.238746 i & -1.001838 + 1.238765 i \\
 40 & (2.6 - 5.8 i)  $\times 10^{37}$ 
 							& -1.001686 + 1.238816 i & -1.001838 + 1.238765 i \\
 50 & (-34. + 3.3 i)  $\times 10^{46}$ 
 							& -1.002011 + 1.238667 i & -1.001838 + 1.238765 i \\
  \ldots & \ldots  & \ldots & \ldots  \\ \hline
    \hline
    \end{tabular}
}
\caption{Resummation of the (divergent) Kapteyn series $\displaystyle\sum_{m\ge 1}\,\dfrac{z^m}m\,J_m(m\epsilon)$ for $z=10\,\exp(\mathrm{i}\pi/3)$.}
\label{Tab.Finale.KS}
\end{table}

\begin{figure}[!ht]
\centering
\begin{minipage}[t]{8cm}
\centerline{\includegraphics[width=8cm]{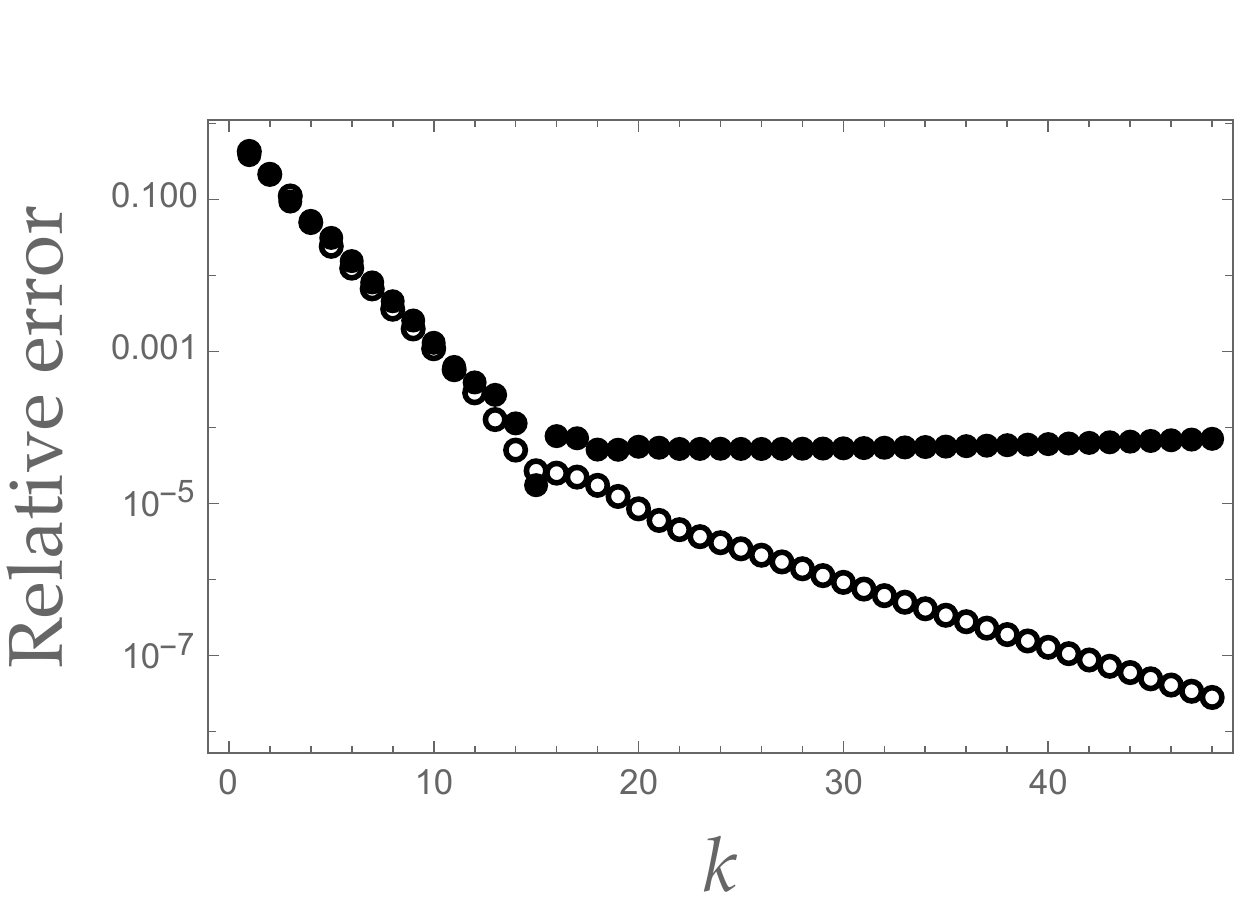}}
\end{minipage}
\caption{Both series in the r.h.s. of Eq.~(\ref{Eq:KE.2.BIS}) have been computed via Levin $d$-  and via Weniger $\delta$-  transformations
of the same order $k$. 
The relative error with respect to the l.h.s. of the same equation,  evaluated through standard method from Eq.~(\ref{Eq:KE.2.TER}),
is plotted against the transformation order for both Levin (black circles) and Weniger (open circles)..}
\label{Fig:Conclusions.1}
\end{figure}

A possible answer can be given on using again Kepler's equation~(\ref{Eq:KE.1}) and, in particular, by suitably extending its ``complexification''
procedure. To this end, the complex parameter $z=\exp(\mathrm{i} M)$ is again introduced, in such a way that $M=-\mathrm{i} \log z$, but 
now $z$ is let to span the whole complex plane. Accordingly, Eq.~(\ref{Eq:KE.3}) takes on the  form
\begin{equation}
\label{Eq:KE.3.BIS}
\begin{array}{l}
\displaystyle
S(\epsilon;M)\,=\,\dfrac 1{\mathrm{i}}\,\sum_{m=1}^\infty\,\frac{z^m\,-\,z^{-m}}m\,J_m(m\,\epsilon)\,,
\end{array}
\end{equation}
and then Eq.~(\ref{Eq:KE.2}) becomes
\begin{equation}
\label{Eq:KE.2.BIS}
\begin{array}{l}
\displaystyle
\sum_{m=1}^\infty\,\frac{z^m}m\,\,J_m(m\,\epsilon)
\,-\,
\sum_{m=1}^\infty\,\frac{z^{-m}}m\,J_m(m\,\epsilon)\,=\,
\Psi\,-\log z\,,
\end{array}
\end{equation}
where now $\Psi$ represents the (complex) solution of a modified version of the Kepler equation, precisely
%
\begin{equation}
\label{Eq:KE.2.TER}
\log z\,=\,\Psi\,-\,\epsilon\,\sinh\,\Psi\,.
\end{equation}
Coming back to the case in Tab.~\ref{Tab.Finale.KS}, in Fig.~\ref{Fig:Conclusions.1}
the result of a numerical experiment carried out on Eqs.~(\ref{Eq:KE.2.BIS}) and~(\ref{Eq:KE.2.TER}) is plotted
for $z=10\,\exp(\mathrm{i}\pi/3)$.
Both series in the r.h.s. of Eq.~(\ref{Eq:KE.2.BIS}) have been computed via Levin $d$-  and via Weniger $\delta$-  transformations
and the relative error with respect to the l.h.s. of the same equation,  evaluated through standard method from Eq.~(\ref{Eq:KE.2.TER}),
is plotted against the transformation order for both Levin (black circles) and Weniger (open circles).

From the figure it appears that the performances of the two transformations are nearly identical for small values of the 
transformation order $k$ (up to about ten). After a relative error of the order of $10^{-5}$ has been reached,  Levin $d$- 
seems to be no longer able to further reducing the error itself, differently from Weniger transformation  which continues to
work, even though with a smaller convergence rate. 
    
Pad\'e approximants, together with Wynn $\epsilon$-algorithm, continue to be the most favorite mathematical tool 
to resum divergent  as well as slowly convergent series in theoretical physics. This is in part due to the fact that a solid
convergence theory for Pad\'e approximants has been developed in the past, especially in the 
case of Stieltjes functions and Stieltjes asymptotic series~\cite{Bender/Orszag/1978,Baker/Graves-Morris/1996,Widder/1938,Widder/1946}. 
Differently from Wynn's $\epsilon$-algorithm, Levin-type transformations
have specifically been designed for resumming alternating factorially divergent asymptotic series.
As a consequence, their numerical performances turn out to be considerably superior with respect to the performances of the former. 
Unfortunately, a convergence theory for Levin-type transformations is still largely missing within the present mathematical
literature, despite their indisputable superiority over Pad\'e in several cases, the probably most 
paradigmatic being the Euler series~\cite{Borghi/2010,Borghi/Weniger/2015}. As a consequence, to build up such a convergence theory 
any  new ``experimental''  result would be helpful. 

In this respect, the deep numerical exploration of the Stieltjes nature of the complex Kapteyn series~(\ref{Eq:KE.4})
operated via  $d$- and $\delta$-transformations should certainly be welcomed, not only to convince skeptic
readers about the extraordinary retrieving capabilities of Levin-type nonlinear sequence transformations. 
Another clue supporting our ``Stieltjness'' conjecture is the exponential convergence displayed by both $d$- and $\delta$-transformations
when they were employed to solve KE via ``brutal force''.

From the perspective of somebody who only wants to employ Levin-type
transformations as computational tools, the situation is not really bad. 
However, our theoretical understanding of the convergence
properties of Levin-type transformations is far from satisfactory, especially with respect to the summation of 
factorially divergent series. In the case of Pad\'{e} approximants, if the input data
are the partial sums of a so-called Stieltjes series, it can be proved rigorously that certain subsequences of the Pad\'{e} 
table converge to the value of the corresponding Stieltjes function~\cite[Chapter 5]{Baker/Graves-Morris/1996}. 
Proving that a divergent power series is a Stieltjes series implies that such a series would be Pad\'{e} summable.
For instance, in~\cite{Bender/Weniger/2001} numerical evidences of the fact that the factorially divergent perturbation expansion for an 
energy eigenvalue of the PT-symmetric Hamiltonian $H(\lambda) =p^2+1/4x^2+\mathrm{i} \lambda x^3$ is a Stieltjes series and thus Pad\'e 
summable, were given. Such a conjecture was later proven rigorously in~\cite{Grecchi/Maioli/Martinez/2009}.

It would be interesting to explore a similar possibility also for Levin-type nonlinear transformations. In other words, under which condtions (if any) 
a Stieltjes series can be summed up to the correct limit by using, for instance, Weniger $\delta$ transformation? 
Presently, it seems not only a definite answer but also weaker conjectures concerning such a question are out of our possibilities.
We hope what is contained in the present paper could stimulate the interest of physicists 
and mathematicians about the importance, both theoretical and practical, of Levin-type nonlinear 
sequence transformations which, in our opinion, have not yet gained the recognition they would definitely deserve. 

\acknowledgments

I wish to thank Michela Redivo-Zaglia and Claude Brezinski for their kind hospitality at the Centre International de Rencontres Math\'ematique
in Luminy, France. I also wish to thank Turi Maria Spinozzi  for his help during the preparation of the manuscript.

\appendix

\section{A recursive algorithm for generating Debye's polynomials}
\label{App:B}

Let us write the explicit expression  of the $k$th-order Debye polynomial $U_k(t)$ as
\begin{equation}
\label{Eq:AppendixB.1}
\begin{array}{l}
\displaystyle
U_k(t)\,=\,\sum_{m=0}^{3k}\,a^k_m\,t^m\,,
\end{array}
\end{equation}
where the superscript $k$ in the symbol $a^k_m$ refers to the polynomial order. Similarly we write
\begin{equation}
\label{Eq:AppendixB.1.1}
\begin{array}{l}
\displaystyle
U_{k+1}(t)\,=\,\sum_{m=0}^{3k+3}\,a^{k+1}_m\,t^m\,,
\end{array}
\end{equation}
so that the task is to find the relationship between the sequences $\{a^k_m\}$ and $\{a^{k+1}_m\}$.
To this aim, all we have to do is to implement the recursive definition~(\ref{Eq:Debye.1.1}) together with 
Eqs.~(\ref{Eq:AppendixB.1}) and~(\ref{Eq:AppendixB.1.1}). First of all we have
\begin{equation}
\label{Eq:AppendixB.2}
\begin{array}{l}
\displaystyle
U'_k(t)\,=\,
\sum_{m=1}^{3k}\,m\,a^k_m\,t^{m-1}\,=\,
\sum_{m=0}^{3k-1}\,(m+1)\,a^k_{m+1}\,t^{m}\,,
\end{array}
\end{equation}
which, after some algebra, leads to
\begin{equation}
\label{Eq:AppendixB.2.1}
\begin{array}{l}
\displaystyle
\dfrac{t^2(1-t^2)}2\,U'_k(t)\,=\,
\sum_{m=2}^{3k+1}\,\dfrac{m-1}2\,a^k_{m-1}\,t^{m}\,\\
\\
\displaystyle
\,-\,
\sum_{m=4}^{3k+3}\,\dfrac{m-3}2\,a^k_{m-3}\,t^{m}\,.
\end{array}
\end{equation}
Similarly  we have
\begin{equation}
\label{Eq:AppendixB.3}
\begin{array}{l}
\displaystyle
\dfrac 18\,\int_0^t\,(1-5\xi^2)\,{U}_k(\xi)\,\d\xi\,=\,\\
\\
\displaystyle
\,=\,
\sum_{m=0}^{3k}\,\dfrac{a^k_m}{8(m+1)}\,t^{m+1}\,-\,
\sum_{m=0}^{3k}\,\dfrac{5a^k_m}{8(m+3)}\,t^{m+3}\,=\,\\
\\
\displaystyle
\,=\,
\sum_{m=1}^{3k+1}\,\dfrac{a^k_{m-1}}{8m}\,t^{m}\,-\,
\sum_{m=3}^{3k+3}\,\dfrac{5a^k_{m-3}}{8m}\,t^{m}\,,
\end{array}
\end{equation}
so that, on taking Eq.~(\ref{Eq:AppendixB.2.1})
into account,
\begin{equation}
\label{Eq:AppendixB.4}
\begin{array}{l}
\displaystyle
\dfrac{t^2(1-t^2)}2\,U'_k(t)\,+\,
\dfrac 18\,\int_0^t\,(1-5\xi^2)\,{U}_k(\xi)\,\d\xi\,=\,\\
\\
\displaystyle
\,=\,
\dfrac{a^k_0}{8}t\,+\,\dfrac{9a^k_1}{16}t^2\,+\,
\dfrac{5}{24}(5a^k_2-a^k_0)t^2\,\\
\\
-\dfrac{3(2k+1)(6k-1)}{8(3k+2)}a^k_{3k-1}\,t^{3k+2}\,\\
\\
-\dfrac{36k(k+1)+5}{24(k+1)}a^k_{3k}\,t^{3k+3}\,\\
\\
\displaystyle
\,+\,
\sum_{m=4}^{3k+1}\,
\left[
\dfrac{(2m-1)^2}{8m}\,a^k_{m-1}\,
\,-\,
\dfrac{4m(m-3)+5}{8m}\,a^k_{m-3}\,
\right]\,t^m\,.
\end{array}
\end{equation}
Finally, on substituting from Eqs.~(\ref{Eq:AppendixB.1.1}) and~(\ref{Eq:AppendixB.4})
into Eq.~(\ref{Eq:Debye.1.1}), after straightforward algebra  the following recurrence
relation is obtained at once:
\begin{equation}
\label{Eq:AppendixB.5}
\begin{array}{l}
\displaystyle
a^{k+1}_0\,=\,0\,\\
\\
\displaystyle
a^{k+1}_1\,=\,\dfrac{a^k_0}{8}\,\\
\\
\displaystyle
a^{k+1}_2\,=\,\dfrac{9a^k_1}{16}\,\\
\\
\displaystyle
a^{k+1}_3\,=\,\dfrac{5}{24}(5a^k_2-a^k_0)\,\\
\\
a^{k+1}_m\,=\,\dfrac{(2m-1)^2}{8m}\,a^k_{m-1}\,
\,-\,
\dfrac{4m(m-3)+5}{8m}\,a^k_{m-3}\,\\
4 \le m \le 3k+1\,\\
\\
\displaystyle
a^{k+1}_{3k+2}\,=\,-\dfrac{3(2k+1)(6k-1)}{8(3k+2)}a^k_{3k-1}\,\\
\\
\displaystyle
a^{k+1}_{3k+3}\,=\,-\dfrac{36k(k+1)+5}{24(k+1)}a^k_{3k}
\end{array}
\end{equation}

\section{Proof of Eq.~(\ref{Eq:Debye.1.1.1.0})}
\label{App:C}

It is sufficient to start from the last row of Eq.~(\ref{Eq:AppendixB.5}), i.e.,
\begin{equation}
\label{Eq:AppendixC.1}
\begin{array}{l}
\displaystyle
a^{k+1}_{3k+3}\,=\,-\dfrac{36k(k+1)+5}{24(k+1)}a^k_{3k}\,,
\end{array}
\end{equation}
which, in the limit of $k \gg 1$, becomes
\begin{equation}
\label{Eq:AppendixC.2}
\begin{array}{l}
\displaystyle
a^{k+1}_{3k+3}\,\sim\,-\dfrac{3}{2}\,k\,a^k_{3k}\,,\qquad\qquad k\to\infty\,,
\end{array}
\end{equation}
and whose explicit solution is
\begin{equation}
\label{Eq:AppendixC.3}
\begin{array}{l}
\displaystyle
a^{k}_{3k}\,\sim\,C\,\left(-\dfrac{3}{2}\right)^{k-1}\,(k-1)!\,,\qquad\qquad k\to\infty\,,
\end{array}
\end{equation}
 i.e., Eq.~(\ref{Eq:Debye.1.1.1.0}).

%
%

%


\begin{thebibliography}{00}

\bibitem{Colwell/1993} P. Colwell,
{\em Solving Kepler's Equation Over Three Centuries}
(Willmann-Bell, Richmond, 1993).



\bibitem{Eisinberg2010215}
A. Eisinberg, G. Fedele, and D. Frascino,
``Kepler's equation and limit cycles in a class of PWM feedback control systems,''
Nonlinear Dynamics, \textbf{62,} 215 - 227  (2010). 

\bibitem{Gu2010716}
S.-L. Gu  and H.-B. Zhang,
``Noether symmetry and the Hojman conserved quantity of the Kepler equation,''
Acta Physica Sinica,  \textbf{59,} 716 - 718  (2010).


\bibitem{Davis201059}
 J.J.  Davis, D. Mortari, and C. Bruccoleri,,
``Sequential solution to Kepler's equation,''
Celestial Mechanics and Dynamical Astronomy, \textbf{108,} 59 - 72  (2010).

\bibitem{Antonov2011182}
V. A. Antonov,  L. N. Sudov,  and K.V. Kholshevnikov,
``Solution of Kepler's equation for rectilinear motion,''
Astronomy Reports, \textbf{55,}  182-186  (2011).

\bibitem{Liu2012}
H.W. Liu,  L.F. Li,  and S.T. Yang,
``Conformal invariance, Mei symmetry and the conserved quantity of the Kepler equation,''
Acta Physica Sinica, \textbf{61,}  200- 202  (2012).

\bibitem{Yang2012517}
L. . Yang  and L. Ma,
``Chaos optimization algorithm for Kepler's equation,''
Journal of University of Shanghai for Science and Technology,
\textbf{34,}  517-519  (2012).

\bibitem{Calvo2013143}
M. Calvo, A.  Elipe, J.I. Montijano, and L. R\'andez,
``Optimal starters for solving the elliptic Kepler's equation,''
Celestial Mechanics and Dynamical Astronomy, \textbf{115,}  143-160  (2013).

\bibitem{Farnocchia201321}
D. Farnocchia,  D.B.  Cioci,  and A. Milani, 
``Robust resolution of Kepler's equation in all eccentricity regimes,''
Celestial Mechanics and Dynamical Astronomy,
\textbf{116,}  21-34  (2013).


\bibitem{Mortari20141}
D. Mortari,  and A. Elipe,
``Solving Kepler's equation using implicit functions,''
Celestial Mechanics and Dynamical Astronomy,
\textbf{118,}  1-11  (2014).

\bibitem{Avendano201427}
M. Avendano, V.  Mart\'an-Molina, and J. Ortigas-Galindo, 
``Solving Kepler's equation via Smale's $\alpha$-theory,''
Celestial Mechanics and Dynamical Astronomy,
\textbf{119,}  27-44  (2014).


\bibitem{Badolati2015316}
E. Badolati, and S. Ciccone, 
``On the history of some explicit formulae for solving Kepler's equation,''
Astronomische Nachrichten,
\textbf{336,}  316-320  (2015).


\bibitem{Lynden-Bell2015363}
D. Lynden-Bell, 
``An approximate analytic inversion of Kepler's equation,''
Monthly Notices of the Royal Astronomical Society,
\textbf{447,}  363-365  (2015).

\bibitem{Avendano2015435}
M. Avendano, M.  V. Mart\'an-Molina, and J. Ortigas-Galindo,
``Approximate solutions of the hyperbolic Kepler equation,''
Celestial Mechanics and Dynamical Astronomy,
\textbf{123,}  435-451  (2015).

\bibitem{Oltrogge2015271}
D.L. Oltrogge, 
``Efficient Solutions of Kepler's Equation via Hybrid and Digital Approaches,''
Journal of the Astronautical Sciences,
\textbf{62,}  271-297  (2015).


\bibitem{Rauh20161}
A. Rauh and J. Parisi,
``Quantum mechanical correction to Kepler's equation,''
Advanced Studies in Theoretical Physics,
\textbf{10,}  1-22  (2016).

\bibitem{Ebaid2016457}
A. Ebaid and A.B. Al-Blowy, 
``Properties of Bessel function solution to Kepler's equation with application to opposition and conjunction of Earth-Mars,''
Zeitschrift fur Naturforschung - Section A Journal of Physical Sciences,
\textbf{71,}  457-464  (2016).

\bibitem{Calvo201719}
M. Calvo, A. Elipe, J.I. Montijano, and L. R\'andez, 
``Convergence of starters for solving Kepler'S equation via Smale's test,''
Celestial Mechanics and Dynamical Astronomy,
\textbf{127,}  19-34  (2017).


\bibitem{Raposo-Pulido20171702}
V. Raposo-Pulido and J. Pel\'aez, 
``An efficient code to solve the Kepler equation. Elliptic case,''
Monthly Notices of the Royal Astronomical Society,
\textbf{467,}  1702-1713  (2017).

\bibitem{Chanda201}
S. Chanda,  G.W. Gibbons, and P. Guha,
``Jacobi-Maupertuis metric and Kepler equation,''
International Journal of Geometric Methods in Modern Physics,
\textbf{14,}  1730002   (2017).

\bibitem{Boetzel/Susobhanan/Gopakumar/Klein/Jetzer/2017}
Y. Boetzel, A. Susobhanan, A. Gopakumar, A. Klein, and P. Jetzer,
``Solving post-Newtonian accurate Kepler equation,''
Phys. Rev. D \textbf{96,} 044011 (2017). 

\bibitem{Ebaid20171}
A. Ebaid, R.  Rach and E. El-Zahar,
``A new analytical solution of the hyperbolic Kepler equation using the Adomian decomposition method,''
Acta Astronautica,
\textbf{138,}  1-9   (2017).

\bibitem{Alshaery2017933}
A. Alshaery, 
``The Homotopy Perturbation Method for Accurate Orbits of the Planets in the Solar System: The Elliptical Kepler Equation,''
Zeitschrift fur Naturforschung - Section A Journal of Physical Sciences,
\textbf{10,}  933-940   (2017).

\bibitem{Alshaery201727}
A. Alshaery and A. Ebaid, 
``Accurate analytical periodic solution of the elliptical Kepler equation using the 
Adomian decomposition method,''
Acta Astronautica,
\textbf{140,}  27-33   (2017).

\bibitem{Elipe2017415}
A. Elipe,  J.I. Montijano, and L. R\'andez, and M. Calvo,
``An analysis of the convergence of Newton iterations for solving elliptic Kepler's equation,''
Celestial Mechanics and Dynamical Astronomy,
\textbf{129,}  415-432   (2017).

\bibitem{Lppez20182583}
R. Lopez, D.  Hautesserres, and  J.F. San-Juan, 
``The solution of the generalized Kepler's equation,''
Monthly Notices of the Royal Astronomical Society,
\textbf{473,}  2583-2589   (2018).

\bibitem{Aljohani2018}
A.F. Aljohani,  R.  Rach, E. El-Zahar, A.M. Wazwaz,  and A. Ebaid, 
``Solution of the hyperbolic kepler equation by adomian's asymptotic decomposition method,''
Romanian Reports in Physics,
\textbf{70,}  14   (2018).

\bibitem{Orlando2018849}
F.G.M. Orlando, C.  Farina, C.A.D. Zarro,  and  P. Terra,
``Kepler's equation and some of its pearls,''
American Journal of Physics,
\textbf{86,}  11   (2018).

\bibitem{Raposo-Pulido2018}
V. Raposo-Pulido. and J. Pel\'aez,
``An efficient code to solve the Kepler equation: Hyperbolic case,''
Astronomy and Astrophysics,
\textbf{619,}  A129   (2018).

\bibitem{Zechmeister2018}
M. Zechmeister,
``CORDIC-like method for solving Kepler's equation,''
Astronomy and Astrophysics,
\textbf{619,}  A128   (2018).

\bibitem{Ibrahim20192269}
R.H. Ibrahim  and A.R.H. Saleh, 
``Re-evaluation solution methods for Kepler's equation of an elliptical orbit,''
Iraqi Journal of Science,
\textbf{60,}  2269-2279   (2019).

\bibitem{Calvo2019}
M. Calvo, A.  Elipe, J.I.  Montijano, and L. R\'andez, 
``A monotonic starter for solving the hyperbolic Kepler equation by Newton method,''
Celestial Mechanics and Dynamical Astronomy,
\textbf{131,}  18   (2019).

\bibitem{Tommasini2020}
D. Tommasini and D.N. Olivieri, 
``Fast switch and spline scheme for accurate inversion of nonlinear functions: The new first 
choice solution to Kepler's equation,''
Applied Mathematics and Computation,
\textbf{364,}  124677  (2020).

\bibitem{Abubekerov/Gostev/2020}
M. K. Abubekerov  and N. Yu. Gostev
``Solution of Kepler's Equation with Machine Precision,''
Astronomy Reports,  \textbf{64,} 1060 - 1066 (2020).

\bibitem{Sacchetti/2020}
A. Sacchetti
``Francesco Carlini: Kepler's equation and the asymptotic solution to singular differential equations,''
Historia Mathematica \textbf{53,} 1- 32 (2020).

\bibitem{Tommasini2020b}
D. Tommasini,  and D.N. Olivieri, 
``Fast Switch and Spline Function Inversion Algorithm
with Multistep Optimization and k-Vector Search for
Solving Kepler's Equation in Celestial Mechanics,''
Mathematics, \textbf{8,}  2017  (2020).

\bibitem{Zechmeister2021}
M. Zechmeister,
``Solving Kepler'ss equation with CORDIC double iterations,''
MNRAS \textbf{500,} 109-117 (2021)

\bibitem{Gonzalez/Hernandez/2021}
O. Gonz\'alez-Gaxiola, S. Hern\'andez-Linares,
``An Efficient Iterative Method for Solving the Elliptical Kepler's Equation,''
Int. J. Appl. Comput. Math.
\textbf{7,} 42(1)- 42(12) (2021).

\bibitem{Tommasini/2021}
D. Tommasini,
``Bivariate Infinite Series Solution of Kepler's Equations,''
Mathematics  \textbf{9,} 785  (2021).

\bibitem{Philcox/Goodman/Slepian/2021}
O. H. E. Philcox , J. Goodman and Z. Slepian,
``Kepler's Goat Herd: An exact solution to Kepler's equation for elliptical orbits,''
Mathematics  \textbf{9,} 785  (2021).

%

\bibitem{Kapteyn/1893} 
W. Kapteyn,
Ann. Sci. de l'Ecole Norm. Sup. Ser. 
\textbf{3,}  91-122 (1893).

\bibitem{Eisinberg/2010}
A. Eisinberg, G. Fedele, A. Ferrise, and D. Frascino, 
``On an integral representation of a class of Kapteyn (Fourier-Bessel) series: Kepler's equation, radiation problems and Meissel's expansion,''
Applied Mathematics Letters, \textbf{23,}  1331-1335 (2010).

\bibitem{Tautz/Lerche/Dominici/2011} 
R. C. Tautz, I. Lerche, and D. Dominici,
``Methods for summing general Kapteyn series,''
J. Phys. A: Math. Theor. \textbf{44} 385202(1)-385202(14) (2011).

\bibitem{Baricz/Jankov/Pogany/2011} 
\'A. Baricz, D. Jankov, and T. K. Pog\'any,
``Integral representation of first kind Kapteyn series,''
J. Math. Phys. \textbf{52,} 043518 (2011).

\bibitem{Nikishov/2014}
A.  I.  Nikishov,
``Kapteyn series and photon emission,''
 Bulletin of the Lebedev Physics Institute, \textbf{41,}  332 - 338 (2014).
 
\bibitem{Dragana/Masirevic/Pogany/2017}
A. Dragana, J. Masirevic, and T. K. Pog\'any,
``Kapteyn Series,'' 
Lecture Notes in Mathematics. \textbf{2207,}  87 - 111 (2017). 

\bibitem{Xue/Li/Man/Xing/Liu/Li/Wu/2019}
X. Xue, Z. Li, Y. Man, S. Xing, Y. Liu, B. Li, and Q. Wu,
``Improved Massive MIMO RZF Precoding Algorithm Based on Truncated Kapteyn Series Expansion,''
Information \textbf{10,} 136 (2019).

\bibitem{Watson/1995} G. Watson,
{\em A Treatise on Theory of Bessel Functions}

\bibitem{Caliceti/Meyer-Hermann/Ribeca/Surzhykov/Jentschura/2007}
E. Caliceti, M. Meyer-Hermann, P. Ribeca,  A. Surzhykov, and U. D. Jentschura,
``From useful algorithms for slowly convergent series to physical
predictions based on divergent perturbative expansions,''
Physics Reports \textbf{446, } 1 - 96,  (2007).

\bibitem{DLMF} {\em NIST Digital Library of Mathematical Functions}. http://dlmf.nist.gov/, Release 1.1.3 of 2021-09-15. F. W. J. Olver, A. B. Olde Daalhuis, D. W. Lozier, B. I. Schneider, R. F. Boisvert, C. W. Clark, B. R. Miller, B. V. Saunders, H. S. Cohl, and M. A. McClain, eds.



\bibitem{Costin/2009}
O. Costin, 
{\em Asymptotics and Borel Summability,} Chapman \& Hall/CRC, Boca Raton (2009).

\bibitem{Baker/Graves-Morris/1996}
G. A. Baker Jr. and P. Graves-Morris, 
{\em Pad\'e Approximants}, 2nd ed., Cambridge U. P., Cambridge, 1996.

\bibitem{Levin/1973}
D. Levin, 
``Development of non-linear transformations for improving convergence of sequences,''
Int. J. Comput. Math. B \textbf{3,} 371-388 (1973).

\bibitem{Wynn/1956}
P. Wynn, 
``On a device for computing the $e_m(S_n)$ transformation,''
Math. Tables Other Aids Comput. \textbf{10,} 91- 96 (1956).

\bibitem{Weniger/1989} E. J. Weniger,
Comput. Phys. Rep. \textbf{10,} 189 - 371 (1989).

\bibitem{Fleckenstein/1941} J.O. Fleckenstein,
``Notiz zur Lagrangeschen Losung des Keplerschen Problems,''
Commentarii Mathematici Helvetici \textbf{13,} 83-89 (1941).

\bibitem{Borghi/Santarsiero/2003}
R. Borghi, M. Santarsiero, 
``Summing Lax series for nonparaxial beam propagation,''
 Opt. Lett. \textbf{28,}  774 - 776 (2003).

\bibitem{borghiOL-07} R. Borghi, 
``Evaluation of diffraction catastrophes by using Weniger transformation,''
Opt. Lett. {\bf 32,} 226 - 228 (2007).

\bibitem{borghiJOSAA-08a}  R. Borghi, 
``{Summing Pauli asymptotic series to solve the wedge problem},''
J. Opt. Soc. Am. A \textbf{25,} 211 - 218 (2008).

\bibitem{borghiJOSAA-08b}  R. Borghi, 
``On the numerical evaluation of cuspoid diffraction catastrophes,''
J. Opt. Soc. Am. A \textbf{25,} 1682 - 1690 (2008).

\bibitem{borghiJOSAA-08c}  R. Borghi, 
``Joint use of the {W}eniger transformation and hyperasymptotics for accurate asymptotic evaluations of a class of saddle-point integrals,''
Phys. Rev. E \textbf{78,} 026703(1) - 026703(11) (2008).

\bibitem{borghiJOSAA-09}  R. Borghi, 
``Joint use of the {W}eniger transformation and hyperasymptotics for accurate asymptotic evaluations of a class of saddle-point integrals.
{II}. {H}igher-order  transformations,''
Phys. Rev. E \textbf{80,} 016704(1) - 016704(15) (2009).

\bibitem{borghiJOSAA-9b}  R. Borghi and M. A. Alonso, 
``Free-space asymptotic far-field series,''
J. Opt. Soc. Am. A \textbf{26,} 2410 - 2417 (2009).

\bibitem{borghiJOSAA-10}  R. Borghi, 
``On the numerical evaluation of umbilic diffraction catastrophes,''
J. Opt. Soc. Am. A \textbf{27,} 1661 - 1670 (2010).

\bibitem{borghiJOSAA-11}  R. Borghi, 
``Numerical evaluation of umbilic diffraction catastrophes of codimension four,''
J. Opt. Soc. Am. A \textbf{28,} 887 - 896 (2011).

\bibitem{borghiOL-11b} R. Borghi, 
``Optimizing diffraction catastrophe evaluation,"
\ol \textbf{36,} 4413 - 4415 (2011).

\bibitem{Borghi/Gori/Guattari/Santarsiero/2011}
R. Borghi, F. Gori, G. Guattari, M. Santarsiero, 
``Decoding divergent series in nonparaxial optics,'' 
Opt. Lett. \textbf{36,} 963 - 965 (2011).

\bibitem{borghiJOSAA-12}  R. Borghi, 
``Numerical computation of diffraction catastrophes with codimension eight,''
Phys. Rev. E \textbf{85,} 046704(1) - 046704(14) (2012).

\bibitem{Borghi/2010} R. Borghi,
``Asymptotic and factorial expansions of {E}uler series truncation errors via exponential polynomials,''
Appl. Numer. Math. \textbf{60,}   1242 - 1250 (2010)

\bibitem{Borghi/Weniger/2015}
``Convergence analysis of the summation of the factorially divergent {E}uler series by {P}ad\'{e}
 approximants and the delta transformation,''
Appl. Numer. Math. \textbf{94,}   149 - 178 (2015)

\bibitem{Bender/Orszag/1978} 
C. M. Bender and  S. A. Orszag, {\em Advanced Mathematical Methods for Scientists and Engineers}, (McGraw-Hill, New York, 1978).

\bibitem{Widder/1938} D. V. Widder, ``The Stieltjes transform,'' 
Trans. Am. Math. Soc. {\bf 43,}  7?60 (1938).

\bibitem{Widder/1946} D. V. Widder, {\em The Laplace Transform,} 
(Princeton U. P., Princeton, 1946). 

\bibitem{Bender/Weniger/2001} C. M. Bender and E. J. Weniger,
``Numerical evidence that the perturbation expansion for a non-Hermitian
PT-symmetric Hamiltonian is Stieltjes,''
J. Math. Phys. \textbf{42,} 2167 (2001).

\bibitem{Grecchi/Maioli/Martinez/2009} 
V. Grecchi, M. Maioli, and A. Martinez, 
``Pad\'e summability of the cubic oscillator,''
J. Phys. A \textbf{42,}  425208 (2009).




\end{thebibliography}
\end{document}